\documentclass[12pt]{article}
\usepackage{a4wide}
\usepackage[english]{babel}
\usepackage{amsfonts}
\usepackage{amssymb, amsthm}
\usepackage{amsmath}
\usepackage[utf8]{inputenc}
\usepackage[all,cmtip]{xy}
\def\beginproof{{\bf Proof: }}
\def\endproof{\hfill$\square$\medskip}
\def\wolp{\stackrel{ o\ } {W_p^1}}
\def\reals{\mathbb R}
\def\lk{\mathcal L_{\mathcal G_k}}

\def\lgk{\mathcal L_{\mathcal G_k}}
\def\lgsk{\mathcal L_{\gsk}}
\def\bgk{\mathfrak B_{\mathcal G_k}}

\def\wlk{{\widetilde{\mathcal L}}_{\gk}}
\def\wlgsk{{\widetilde{\mathcal L}}_{\gsk}}
\def\myliminf{\mathop{\underline{\rm lim}}}
\def\mylimsup{\mathop{\overline{\rm lim}}}
\def\supp{\mathop{\rm supp}}
\def\mes{{\rm mes}}
\def\gk{\mathcal G_k}
\def\gsk{\mathcal G'_k}
\def\wgk{\widetilde{\mathcal G_k}}
\def\wgsk{\widetilde{\mathcal G'_k}}
\def\bwgk{\mathfrak B_{\wgk}}
\def\lwgk{\mathcal L_{\wgk}}
\def\wlwgk{\widetilde{\mathcal L_{\wgk}}}
\def\lwgsk{\mathcal L_{\wgsk}}
\def\wlwgsk{\widetilde{\mathcal L_{\wgsk}}}
\def\vec#1{\left ( \begin{array}{c} #1 \end{array} \right )}

\newtheorem{thm}{Theorem}[section]
\newtheorem{lemma}{Lemma}[section]
\newtheorem{corollary}{Corollary}[section]
\newtheorem{remark}{Remark}[section]

\newtheorem{proposition}{Proposition}[section]

\begin{document}

\title{\normalsize
\bf Multiplicity of 1D-concentrating positive solutions\\ to the Dirichlet problem for equation with $p$-Laplacian}
\author { {\it S.B.~Kolonitskii}\thanks {Supported by RFBR grant 11-01-00825.},\\
Saint-Petersburg State University, \\
{\small e-mail:\ sergey.kolonitskii@gmail.com}}
\date{2012}
\maketitle

\section*{Introduction}

The qualitative theory of quasilinear elliptic equations with variational structure has been a rapidly developing field for several last decades.
In this paper we study the multiplicity effect for a model boundary value problem
\begin{equation}
 -\Delta_p u = u^{q-1} \quad\mbox{in}\quad{\Omega_R};\qquad u>0 \quad\mbox{in}\quad {\Omega_R};\qquad u=0 \quad\mbox{on}\quad \partial{\Omega_R},\label{p_laplacian_eq}
\end{equation}
where ${\Omega_R} = B^n_{R+1} \setminus B^n_{R-1}$ is a spherical layer in $\reals^n$, $\Delta_p u = \mbox{div} \left ( |\nabla u|^{p-2} \nabla u \right )$ is $p$-Laplacian, $1<p<\infty$, $q>p$ and $R>2$. Denote by $p^*_n$ the critical Sobolev exponent, which can be defined from equation $\frac 1 {p^*_n} = \left ( \frac 1 p - \frac 1 n \right )_+$.

Multiplicity of solutions to the problem \eqref{p_laplacian_eq} in the expanding annuli was studied in a number of papers, starting with a seminal paper \cite{coffman}. The term ``multiplicity'' means that for any natural $K$ there exists $R_0 = R_0(p,q,K)$ such that for all $R > R_0$ the problem (\ref{p_laplacian_eq}) has at least $K$ nonequivalent (i.e. not obtainable from each other by rotations) solutions.  Multiplicity in the case $p=2$ was proven in \cite{coffman} for $n=2$ and in \cite{li} for $n \geqslant 4$. In these articles solutions of the boundary problem (\ref{p_laplacian_eq}) were obtained as points of global minimum of functional
\begin{equation}
 J[u] = {\int\limits_{\Omega_R} |\nabla u|^p dx} \Bigg / {\Big ( \int\limits_{\Omega_R} |u|^q dx \Big )^{\frac p q}}\label{J_functional}
\end{equation}
on spaces of functions with certain symmetries. The case $n=3$ was treated in \cite{byeon}. Here, the solutions are obtained as local minimizers of functional \eqref{J_functional}. Solutions obtained in these papers concentrate in the neighbourhood of discrete set of points. Solutions concentrating in the neighboourhood of manifolds were constructed in \cite{Malchiodi}. Many other papers are devoted to semilinear case, but their methods cannot be applied directly to quasilinear problem \eqref{p_laplacian_eq} in case of general $p \in (1,\infty)$.

Multiplicity of solutions for equation with $p$-Laplacian was proven in \cite{ain1} for $n=2$ and in \cite{ain2} for $n \geqslant 4$. Three-dimensional case was treated in \cite{kol2}. By solution hereafter we mean a weak solution of class $\mbox{$\wolp({\Omega_R})$}$. As in semilinear case, solutions constructed in these papers concentrate in the neighbourhood of discrete set of points. Note that multiplicity of solutions constructed in all mentioned papers arises only for $q \in (p, p_n^*)$. The existence of several nonequivalent solutions concentrating in the neighbourhood of manifolds was proven in \cite{ain2}.

We establish the multiplicity of solutions to the problem (\ref{p_laplacian_eq}) in case of $n = 4$ or $n \geqslant 6$, arbitrary $p$ and $ q \in (p, p^*_{n-1})$ (note that $p_{n-1}^* > p_n^*$ if $p<n$). Solutions we will present concentrate in the neighbourhood of certain curves.
\medskip

The structure of paper is as follows: Section 1 contains some auxilary lemmas and definitions. Multiplicity of solutions in $\reals^4$ is obtained in Sections 2,3 and 4. Sections 5 and 6 cover the multiplicity effect in $\reals^n$, $n \geqslant 6$. A more general construction of multiple solutions is considered in Section 7.

Let us introduce some notation. $f = o_R(1)$ means that $\lim\limits_{R\to\infty} f = 0$. All functions subscripted with $R$ have a support inside $\overline{\Omega_R}$ and are continued to $\reals^n \setminus {\Omega_R}$ with zero values. All sequences of radii $R_j$ are considered increasing with respect to $j$ and tending to infinity as $j \to \infty$. When we consider sequences $u_{R_j}$ or $x_{R_j}$, we omit the index $j$.

The ball of radius $\rho$ with center $x$ in $\reals^n$ will be denoted $B_\rho^n(x)$. If $x=0$, we will omit the center and if dimension of containing space is evident from context, we will omit the dimension.

We use the standard notation $L_p(\Omega)$, $W_p^1(\Omega)$, $\wolp(\Omega)$ for the functional spaces. $\|\cdot\|_p$ stands for the norm in $L_p(\reals^n)$.

$L_q(\omega, \mathfrak w)$ is the weighted Lebesgue space, i.e. the space of measurable functions such that \break $\int\limits_\omega |u|^q \mathfrak w dx < \infty$. $W_p^1(\omega,\mathfrak w)$ is the weighted Sobolev space, i.e. the space of weakly differentiable functions such that $\int\limits_\omega (|\nabla u|^p + |u|^p) \mathfrak w dx < \infty$.

Let $\mathcal G$ be a closed subgroup of $\mathcal O(n)$. We define by $\mathfrak B_{\mathcal G} (x, \rho)$ a $\mathcal G$-tube:
$$
 \mathfrak B_{\mathcal G}(x, \rho) = \left\lbrace gy \big | g \in \mathcal G, |x-y| < \rho \right\rbrace.
$$
$\mathcal L_{\mathcal G}$ stands for the space of $\mathcal G$-invariant functions from $\wolp(\Omega_R)$.

Consider a $1$-homogeneous map $P: \Omega_R \to \Omega_R' = P \Omega_R$. Points in $P \Omega_R$ will be denoted as $P x$. Given $\mathcal G \subset \mathcal O(n)$, suppose that $P x = P y$ only if $x$ and $y$ are equivalent up to action of $\mathcal G$, i.e. there exists $g \in \mathcal G$ such that $y = gx$. Consider a $\mathcal G$-invariant function $f: \Omega_R \to \reals$. The function $Pf:\Omega_R' \to \reals$ defined by $(Pf)(Px) = f(x)$ will be called {\bf transplant of function $f$ under action of map $P$}. 
As $P$ is homogeneous, $P \reals^n$ is a cone and $P \Omega_R$ is a part of spherical layer in appropriate euclidean space.

$I_m$ is the identity $m \times m$-matrix and $T_\varphi$ is the $2 \times 2$-rotation matrix:
$$
 T_\varphi = \left ( \begin{array}{cc}
          \cos(\varphi) & \sin(\varphi) \\
          -\sin(\varphi) & \cos(\varphi)
         \end{array}
 \right ).
$$

Constants depending only on $n$, $p$ and $q$ are denoted by $C$. 

\section{Auxilary lemmas and definitions}

This lemma is a variant of Lemma 1.1 in \cite[part 1]{Lio}.
\begin{lemma}
 \label{point_concentration_lemma}
 Let $G(s)$ be a positive function. Consider a sequence $u_j$ of $\mathcal G$-invariant functions and suppose that $\int_{{\reals^n}} G(u_j) dx$ is finite for all $j$. Consider a sequence of points in $x_j \in \reals^n$. Then, passing to subsequence if necessary, one of two following statements holds:

1) ({\it concentration}) There exists $\lambda\in (0,1]$ such that for any $\varepsilon > 0$ there exist $\rho>0$ and $j_0$ such that for any $j \geqslant j_0$ there exists a sequence $\rho'(j)$ tending to infinity such that 
\begin{multline}
 \label{concentration}
 \left | \int\limits_{\mathfrak B_{\mathcal G}(x_j,\rho)} G(u_j) dx - \lambda \int\limits_{{\reals^n}} G(u_j) dx \right | + \\
 + \left | \int\limits_{\reals^n \setminus \mathfrak B_{\mathcal G}(x_j,\rho'(j))} \! G(u_j) dx - (1 - \lambda) \int\limits_{{\reals^n}} G(u_j) dx \right | < \varepsilon \int\limits_{{\reals^n}} G(u_j) dx.
\end{multline}

2) ({\it vanishing}) For all $\rho>0$ the following statement holds:
\begin{equation}
 \label{vanishing}
 \lim_{j\to\infty} \frac {\int\limits_{\mathfrak B_{\mathcal G}(x_j,\rho)} G(u_j) dx} {\int\limits_{\reals^n} G(u_j) dx} =0.
\end{equation}
\end{lemma}

\beginproof
 Let $a_j = \int\limits_{\reals^n} G(u_j) dx$. Consider functions
 \begin{equation*}
  h(\rho,j) = \frac 1 {a_j} \int\limits_{\mathfrak B_{\mathcal G}(x_j,\rho)} G(u_j) dx \qquad\mbox{and}\qquad i(\rho) = \mylimsup_{j\to\infty} h(\rho,j).
 \end{equation*}
 $i(\rho)$ is a nondecreasing and bounded function of $\rho$, so we can consider $\lambda = \lim\limits_{\rho\to \infty} i(\rho)$. If $\lambda=0$, then for all $\rho>0$ $i(\rho)=0$ holds, i.e. statement (\ref{vanishing}) (vanishing) holds.

 Suppose that $\lambda \in (0,1]$. Fix $\varepsilon > 0$. There exists such $\rho_0$ that $i(\rho) \geqslant \lambda-\frac \varepsilon 6$, if $\rho > \rho_0$. Fix such $\rho > \rho_0$.  Then by $\overline{\lim\limits_{j\to\infty}} h(\rho,j) \geqslant \lambda - \frac \varepsilon 4$ we obtain that there exists a subsequence $j_m$ such that $j_m \to \infty$ and $h(\rho, j_m) \geqslant \lambda - \frac {2 \varepsilon} 4$, i.e. 
 \begin{equation}
  \label{loc_dichotomy_1}
  \frac 1 {a_j} \int\limits_{\mathfrak B_{\mathcal G}(x_j,\rho)} G(u_j) dx \geqslant \lambda - \frac {2 \varepsilon} 4.
 \end{equation}
 On the other hand, $i(\rho) \leqslant \lambda$ for all $\rho$. Let sequence $\rho_m$ tend to infinity. Then for any $m$ there exists $j_{0,m}$ such that $h(\rho_m,j) \leqslant \lambda + \frac \varepsilon 2$ for all $j>j_{0,m}$, i.e.
 \begin{equation}
  \label{loc_dichotomy_2}
  \frac 1 {a_j} \int\limits_{B_\mathcal G(x_{j_m},\rho_m)} G(u_{j_m}) dx \leqslant \lambda + \frac \varepsilon 2.
 \end{equation}
 Thus
 \begin{equation*}
  \left | \frac 1 {a_j} \int\limits_{B_\mathcal G(x_{j_m},\rho)} G(u_{j_m}) dx - \lambda \right | \leqslant \frac \varepsilon 2.
 \end{equation*}
 From (\ref{loc_dichotomy_2}) we obtain $\frac 1 {a_j} \int\limits_{\reals^n \setminus B_\mathcal G(x_{j_m},\rho_m)} G(u_{j_m}) dx \geqslant 1 - \lambda - \frac \varepsilon 2$. On the other hand, from (\ref{loc_dichotomy_1}) we obtain $\frac 1 {a_j} \int\limits_{\reals^n \setminus B_\mathcal G(x_{j_m},\rho_m)} G(u_{j_m}) dx \leqslant \frac 1 {a_j} \int\limits_{\reals^n \setminus B_\mathcal G(x_{j_m},\rho)} G(u_{j_m}) \leqslant 1 - \lambda + \frac \varepsilon 2$. Thus
 \begin{equation*}
  \left | 1 - \lambda - \frac 1 {a_j} \int\limits_{\reals^n \setminus B_\mathcal G(x_{j_m},\rho_m)} G(u_{j_m}) dx \right | \leqslant \frac \varepsilon 2.
 \end{equation*}
\endproof

\begin{remark}
 Sequence $x_j$ for which inequality \eqref{concentration} holds is called {\bf concentration sequence for sequence $u_j$}.
\end{remark}

For brevity, we will say ``sequence $u_j$ has concentration sequence $x_j$'' instead of ``sequence $u_j$ has subsequence $u_{j_k}$ such that $x_{j_k}$ is concentration sequence for $u_{j_k}$''.


\begin{remark}
 \label{slow_exterior_radius_remark}
 Sequence of radii $\rho'_j$ can be replaced with any sequence of radii $\rho'_j$ for which $\rho'_j \leqslant \rho_j$ for all $j$ and $\rho'_j \to \infty$ as $j \to \infty$ and $\rho$ can be replaced with any greater constant $\rho'$ without breaking inequality (\ref{concentration}). Obviously we can assume without loss of generality that $\rho'_j$ tends to infinity as slowly as desired.
\end{remark}


\begin{remark}
 \label{equivalent_sequence_remark}
 Two sequences of points $x_j$ and $y_j$ are called {\bf equivalent} if there exists $C>0$ such that for all $j$ ${\rm dist}(\mathcal G x_j,\mathcal G y_j) = \min_{g \in \mathcal G} |x_j - g y_j| < C$.
 
 Increasing $\rho$ and decreasing $\rho'_j$ if needed we can ensure that application of Lemma \ref{point_concentration_lemma} to equivalent sequences $x_j$ and $y_j$ yields either (\ref{vanishing}) or (\ref{concentration}) with the same $\lambda$ simultaneously. Hence in what follows we will consider concentration sequences up to equivalence.
\end{remark}

\begin{lemma}
 \label{vanishing_radius_lemma}
 Let $u_j$ be a sequence of $\mathcal G$-invariant functions. Suppose that $x_j$ is not a concentration sequence for $u_j$. Fix $\varepsilon>0$. Then there exists a sequence of radii $\rho_j \to \infty$ as $j \to \infty$ such that 
 \begin{equation}
  \label{vanishing_varepsilon}
  \int\limits_{\mathfrak B_{\mathcal G}(x_j,\rho_j)} G(u_j) dx < \varepsilon \int\limits_{\reals^n} G(u_j) dx.
 \end{equation}
\end{lemma}

\beginproof
 Define
 \begin{equation*}
  F(j,r)= \frac {\int\limits_{\mathfrak B_{\mathcal G}(x_j,r)} G(u_j) dx} {\int\limits_{\reals^n} G(u_j) dx}.
 \end{equation*}
 Take $j_1$ such that for any $j>j_1$ inequality $F(j,1) < \varepsilon$ holds. Select $r_0>0$ such that $F(r_0) < \varepsilon$ for all $j \leqslant j_0$ inequality $F(j,r_0) < \varepsilon$ holds and put $\rho_j = r_0$ for $j \leqslant j_1$.

 Now take $j_2$ such that for any $j>j_2$ inequality $F(j,2) < \varepsilon$ holds. Without loss of generality we can assume that $j_2 > j_1$. Put $\rho_j = 1$ for $j \in (j_1, j_2]$.

 Suppose now that we have constructed $\rho_j$, $j \leqslant j_k$ such that 1) $F(j, \rho_j) < \varepsilon$ for any $j \leqslant j_k$, 2) $F(j,k) < \varepsilon$ for any $j > j_k$, 3) $j_k > k$, and 4) $\rho_{j_k} = k-1$. Take $j_{k+1} > j_k$ such that inequality $F(j,k+1) < \varepsilon$ holds for any $j > j_{k+1}$ and put $\rho_j = k$ for $j \in (j_k, j_{k+1}]$. By induction, we obtain a sequence $\rho_j$ such that \eqref{vanishing_varepsilon} holds.
\endproof

From now on we will assume that all sequences $x_j$ are tending to infinity and their projections on unit sphere converge. By convention, if sequence is indexed with $R_j$ rather then $j$, we assume that $|x_{R_j}| = R_j$. We recall that by convention $R_j \to \infty$ as $j \to \infty$.

Next two statements are well known and we state them without proof.

\begin{proposition}
 (Lemma 1.4 \cite{ain2})
 \label{Friedrichs_Nazarov_proposition}
 Let  $\mathcal D = \{ (R,\Theta)| R \in (R-1,R+1), \Theta \in S \}$ be the intersection of spherical layer with a cone built upon the set $S \subset \mathbb S_n$. Let $p < q < p^*$, $u \in \wolp(\Omega)$, $v = u|_{\mathcal D}$. Assume that $S$ is  ``large enough'' in following sense: any straight line parallel to coordinate axis and intersecting $\mathcal D$ intersects $\partial \mathcal D$ at $\{R-1,R+1\}\times S$, i.e. where $v=0$. Then
\begin{equation*}
 \|v\|_q \leqslant c_0 \|\nabla v\|_p,
\end{equation*}
where $c_0$ does not depend on $R$. If $p \neq n$, this inequality also holds for $q=p^*$.
\end{proposition}
\begin{remark}
 Any domain $S \subset \mathbb S_n$ satisfies conditions of Proposition \ref{Friedrichs_Nazarov_proposition} if $R$ is large enough.
\end{remark}

\begin{proposition}
 (Lemma 3 \cite{kol2})
 \label{vanishing_proposition}
Let sequence $u_R$ be bounded in $\wolp(\Omega_R)$, let $q \in (p,p^*)$, and let for some $\rho>0$
\begin{equation*}
 \lim_{R\to\infty} \sup_{x \in \omega }  {\int\limits_{B(x,\rho)} |u_R|^q dx}  = 0,
\end{equation*}
where $\omega $ is an open set in $\reals^n$, such that $|\omega \cap B(x_0,\rho)| \geqslant A \rho^n$ for all $x \in \omega$ with constant  $A > 0$ independent of $\rho$ and $x_0$. Then $\int\limits_\omega |u_R|^q dx \to 0$ as $R \to \infty$.
\end{proposition}


Next two lemmas allow to freeze weight functions in integral functionals. First lemma is elementary, so we give it without proof.

\begin{lemma}
 \label{weight_freeze_lemma}
 Let $v_j \in L_1(\reals^n,\mathfrak w)$. Suppose that for some $x_j \in \supp v_j$ following condition holds:
 \begin{equation*}
  \frac {\max\limits_{x \in \supp v_j} |\mathfrak w(x) - \mathfrak w(x_j)|} {\min\limits_{x \in \supp v_j} \mathfrak w(x)} \to 0 \quad \mbox{as} \quad j \to \infty.
 \end{equation*}
 Then 
 \begin{equation*}
  \int\limits_{\reals^n} v_j \mathfrak w dx = \mathfrak w(x_j) \int\limits_{\reals^n} v_j dx \cdot (1 + o_j(1)).
 \end{equation*}
\end{lemma}

\begin{lemma}
 \label{matrix_freeze_lemma}
 Let $p >1$. Let $\mathfrak w: \mathcal D \subset \reals^n \to \reals$ be a nonnegative function. Let $A,B_j : \mathcal D \to \reals^{n \times n}$ be nonnegative matrix functions. Let $v_j$ be a sequence of functions such that
 $$
  \int_{\supp v_j} \left ( A(x) \nabla v_j, \nabla v_j \right )^{\frac p 2} \mathfrak w(x) dx < \infty.
 $$
Assume that 
 $$
  |((B_j(x)-A(x))\zeta,\zeta)| \leqslant (A(x)\zeta,\zeta) \cdot o_j(1)
 $$
 for all $x \in \supp v_j$ and $\zeta \in \reals^n$ with some $\alpha_j \to 0$ as $j \to \infty$. 

 Then 
 \begin{equation*}
  \int\limits_{\supp v_j} \left (B_j (x) \nabla v_j,\nabla v_j \right )^{\frac p 2} \mathfrak w(x) dx =  \int\limits_{\supp v_j} \left (A (x) \nabla v_j,\nabla v_j \right )^{\frac p 2} \mathfrak w(x) dx \cdot (1 + o_j(1)).
 \end{equation*}
\end{lemma}

\beginproof
 We will need the following inequality: 
\begin{equation}
 \label{vector_inequality}
 |a+b|^s \leqslant |a|^s + C(s) \left ( |a|^{s-1} |b| + |b|^s \right ).
\end{equation}
 It holds for all  $a,b \in \reals^N$, $s>1$. It follows that 
 \begin{equation}
  \left | |b|^s - |a|^s \right | \leqslant C(s) \left ( |a|^{s-1} |b-a| + |b-a|^s \right )
 \end{equation}

 Suppose that $p \geqslant 2$. Inequality \eqref{vector_inequality} provides that 
 \begin{multline*}
  \left | \int\limits_{\supp v_j} \left (B_j (x) \nabla v_j,\nabla v_j \right )^{\frac p 2} \mathfrak w(x) dx -  \int\limits_{\supp v_j} \left (A (x) \nabla v_j,\nabla v_j \right )^{\frac p 2} \mathfrak w(x) dx  \right | \leqslant \\
   C \Bigg ( \int\limits_{\supp v_j} |((B_j(x) - A(x))\nabla v_j,\nabla v_j)| \left ( A(x) \nabla v_j, \nabla v_j \right )^{\frac p 2 - 1} \mathfrak w(x) dx + \\
  +\int\limits_{\supp v_j} |((B_j(x) - A(x))\nabla v_j,\nabla v_j)|^{\frac p 2} \mathfrak w(x) dx \Bigg )
 \end{multline*}
 Consider the last term:
 \begin{equation*}
  \int\limits_{\supp v_j} |((B_j(x) - A(x))\nabla v_j,\nabla v_j)|^{\frac p 2} \mathfrak w(x) dx \leqslant 
   \int\limits_{\supp v_j} (A(x)\nabla v_j,\nabla v_j)^{\frac p 2} \mathfrak w(x) dx \cdot o_j(1).
 \end{equation*}
 Applying H\"older inequality to first term, we obtain
 \begin{multline*}
  \int\limits_{\supp v_j}  |((B_j(x) - A(x))\nabla v_j,\nabla v_j)| \left ( A(x) \nabla v_j, \nabla v_j \right )^{\frac p 2 - 1} \mathfrak w(x) dx \leqslant \\
  \left ( \int\limits_{\supp v_j} |((B_j(x) - A(x))\nabla v_j,\nabla v_j)|^{\frac p 2} \mathfrak w(x) dx \right )^{\frac 2 p} \left ( \int\limits_{\supp v_j} (A(x)\nabla v_j,\nabla v_j)^{\frac p 2} \mathfrak w(x) dx \right )^{1 - \frac 2 p} \leqslant \\
  \leqslant \int\limits_{\supp v_j} (A(x)\nabla v_j,\nabla v_j)^{\frac p 2} \mathfrak w(x) dx \cdot o_j(1).
 \end{multline*}
 Estimate is proven for $p \geqslant 2$. If $1<p<2$, we use concavity of function $f(t) = t^{\frac p 2}$.
 \begin{multline*}
  \left | \int\limits_{\supp v_j} \left (B_j (x) \nabla v_j,\nabla v_j \right )^{\frac p 2} \mathfrak w(x) dx -  \int\limits_{\supp v_j} \left (A (x) \nabla v_j,\nabla v_j \right )^{\frac p 2} \mathfrak w(x) dx \right | \leqslant \\
  \leqslant {\frac p 2} \int\limits_{\supp v_j}  |((B_j(x) - A(x))\nabla v_j,\nabla v_j)| \left ( A(x) \nabla v_j, \nabla v_j \right )^{\frac p 2 - 1} \mathfrak w(x) dx \leqslant \\
  \leqslant \int\limits_{\supp v_j} (A(x)\nabla v_j,\nabla v_j)^{\frac p 2} \mathfrak w(x) dx \cdot o_j(1).
 \end{multline*}
\endproof

\begin{remark}
 \label{weight_freeze_remark}
 Let $x_j$ be a sequence of points such that $|x_j| \to \infty$ as $j \to \infty$. Let $\rho_j$ be such a sequence that $\frac {\rho_j} {|x_j|} \to 0$ as $j \to \infty$. Let $v_j$ be a sequence of functions such that $\supp v_j \subset B_{x_j,\rho_j}$.

 1) Let $\mathfrak w$ be a positive $m$-homogeneous function. Then conditions of Lemma \ref{weight_freeze_lemma} are satisfied.

 2) Let $A, B_j$ be positive $0$-homogeneous matrix functions such that $A(x_j) = B_j(x_j)$. Then conditions of Lemma \ref{matrix_freeze_lemma} are satisfied.
\end{remark}

\begin{lemma}
 \label{OR_lemma}
 Let $0 < \mathcal A \leqslant \mathcal B < \mathcal C$ and $0<r<1$. Let $\delta = \frac {\mathcal A}  {\mathcal A + \mathcal B}$. Consider a following minimization problem: $f(\alpha,\beta,\gamma) = \mathcal A \alpha^r + \mathcal B \beta^r + \mathcal C \gamma^r \to \min$ with constraints
 \begin{gather}
  1 - \varepsilon \leqslant \mathcal A \alpha + \mathcal B \beta + \mathcal C \gamma \leqslant 1; \nonumber\\
  \mathcal A \alpha \leqslant \delta;\qquad \mathcal A \alpha + \mathcal C \gamma \geqslant \delta - \varepsilon; \label{OR_l_constraints}\\
  \alpha \geqslant 0; \beta \geqslant 0; \gamma \geqslant 0. \nonumber
 \end{gather}
 Then there exists such $\varepsilon_0>0$ depending only on $\mathcal A$, $\mathcal B$, $\mathcal C$ and $r$ that for any  $\varepsilon \in (0, \varepsilon_0)$ minimal value of function $f$ with constraints (\ref{OR_l_constraints}) is strictly greater than $\mathcal B^{1-r}$.

 The same statement holds if $0 < \mathcal A \leqslant \mathcal B$ and $\mathcal C = 0$.
\end{lemma}
\beginproof
 Consider an auxilary problem: minimize function $f$ with constraints
 \begin{gather}
  \mathcal A \alpha + \mathcal B \beta + \mathcal C \gamma = 1; \nonumber\\
  \mathcal A \alpha \leqslant \mathcal \delta;\qquad \mathcal A \alpha + \mathcal C \gamma \geqslant \delta; \label{OR_la_constraints}\\
  \alpha \geqslant 0; \beta \geqslant 0; \gamma \geqslant 0. \nonumber
 \end{gather}
 If $\mathcal C \not = 0$, the domain given by constraints (\ref{OR_la_constraints}) is a flat convex quadrangle.
 $f$ is a concave function, hence minimum is attained at one of corners. Values of function $f$ at corners are
 \begin{align*}
  f_1 &= \mathcal C^{1-r};\\
  f_2 &= (1 - \delta)^r \mathcal B^{1-r} + \delta^r \mathcal C^{1-r};\\
  f_3 &= (1 - \delta)^r \mathcal C^{1-r} + \delta^r \mathcal A^{1-r};\\
  f_4 &= (1 - \delta)^r \mathcal B^{1-r} + \delta^r \mathcal A^{1-r}.
 \end{align*}
Obviously $f_4 < f_2$ and $f_4 < f_3$. $f_4$ as a function of $\delta$ attains maximum when $\delta = \frac {\mathcal A}  {\mathcal A + \mathcal B}$, and its maximal value is $(\mathcal A + \mathcal B)^{1-r}$. Then the minimum of function $f$ under constraints (\ref{OR_la_constraints}) is equal to $\mathfrak M :=  \min \left ( (\mathcal A + \mathcal B)^{1-r}, \mathcal C^{1-r} \right ) > \mathcal B^{1-r}$. 

If $\mathcal C=0$, constraints \eqref{OR_la_constraints} imply that $\alpha = \frac \delta {\mathcal A}$, $\beta = \frac {1 - \delta} {\mathcal B}$, i.e. domain defined by these constraints consists of exactly one point. Value of function $f$ at this point is $f_4$, which is greater than $\mathcal B^{1 - r}$.

Constraints (\ref{OR_l_constraints}) are continuous with respect to $\varepsilon$, so the domain given by these constraints is close to the domain given by (\ref{OR_la_constraints}), and $f$ is continuous, hence its minimum under constraints (\ref{OR_l_constraints}) is close to $\mathfrak M$ if $\varepsilon$ is small enough.
\endproof

\begin{lemma} 
 \label{hump_destruction_lemma}
 Let $1<p<q<\infty$. Let functions $a, b, c \in \wolp(\reals^n) \cap L_q(\reals^n)$ have separated supports. Suppose that $b \neq 0$, $c \neq 0$, and 
 \begin{equation}
  \label{hump_positive_derivative_condition}
  \|\nabla b\|_p^p/\|b\|_q^q \geqslant \|\nabla c\|_p^p/\|c\|_q^q.
 \end{equation}
 Suppose that function $u = a+b+c$ is normalized in $L_q(\reals^n)$. Then there exists such function $U$ that 
\begin{gather*}
 U=a \mbox{ if } x \in \reals^n\setminus \supp c;\\
 \int \limits_{\reals^n} U^q dx = \int \limits_{\reals^n} u^q dx = 1;\\
 \|\nabla U\|_p^p < \|\nabla u\|_p^p - C(b,c); 
\end{gather*}
\end{lemma}

\beginproof
Let 
\begin{equation*}
 t_0  = \frac {\|b\|_q^q} {\|b\|_q^q + \|c\|_q^q}.
\end{equation*}
Consider a function $u_t = a + \left ( \frac t {t_0} \right )^{\frac 1 q} b + \left ( \frac {1-t} {1-t_0} \right )^{\frac 1 q} c$. Then $u_{t_0} = u$ and, since supports of summands are separated from one another,
\begin{math}
 \|u_t\|_q^q = 1.
\end{math}

 Consider the function
 \begin{equation*}
  f(t) = \|\nabla u_t\|_p^p = \int\limits_{\reals^n} |\nabla a|^p dx + \left ( \frac t {t_0} \right )^{\frac p q} \int\limits_{\reals^n} |\nabla b|^p dx + \left ( \frac {1-t} {1-t_0} \right )^{\frac p q} \int\limits_{\reals^n}  |\nabla c|^p dx.
 \end{equation*}
 Then
 \begin{equation*}
  f'(t_0) = \frac p q  \left [ \frac 1 {t_0} \|\nabla b\|_p^p - \frac 1 {1 - t_0} \|\nabla c\|_p^p \right ] = \frac p q (\|b\|_q^q + \|c\|_q^q)  \left [ \|\nabla b\|_p^p/\|b\|_q^q - \|\nabla c\|_p^p/\|c\|_q^q \right ] \geqslant 0.
 \end{equation*}
 Since $p<q$, $f$ is obviously strictly concave. Therefore it is strictly increasing on $[0,t_0]$. Thus $\|\nabla u_0\|_p^p = f(0) < f(t_0) - C(b,c) = \|\nabla u\|_p^p - C(b,c)$. 

 A following bound can be obtained after certain evaluations:
 \begin{equation*}
 f(0) < f(\frac {t_0} 4) \leqslant f(t_0) - \frac p {16 q} \left ( 1 - \frac p q \right ) \left [ \|\nabla b\|_p^p + \left ( \frac {2 \|c\|_q^q} {2 \|c\|_q^q + \|b\|_q^q} \right )^{2 -\frac p q} \frac {\|b\|_q^{2q}} {\|c\|_q^{2q}} \|\nabla b\|_p^p \right ].
 \end{equation*}
\endproof

\section{4D case. Global and local profiles}

In this section we consider the case $n=4$. Let group $\gk$, $k \geqslant 2$ be generated by matrices
\begin{equation}
 \mathcal T = \left (
  \begin{array}{cc}
    0 & I \\
    I & 0 \\
  \end{array}
 \right ),
 \mathcal S_k = \left (
  \begin{array}{cc}
    T_{\frac {2 \pi} k} & 0 \\
    0 & I\\
  \end{array}
 \right )
 \mbox{ and }
 \mathcal R_\varphi = \left (
  \begin{array}{cc}
    T_\varphi & 0 \\
    0 & T_\varphi\\
  \end{array}
 \right ),
\mbox{ where } \varphi \in [0,2 \pi].
\end{equation}

Let $\mathcal G_0$ be the group of all $\mathcal R_\varphi$. Let $\mathcal H_k$ be the group generated by matrices $\mathcal T$ and $\mathcal S_k$.

Orbits of points under action of group $\gk$ have dimension $1$, but their topological structure may differ. There are three classes of points:
\begin{enumerate}
 \item Points that have coordinates $(x_1,x_2,0,0)$ or $(0,0,x_3,x_4)$. Orbits of such points have length $4 \pi |x|$. Their orbits will be called {\bf polar} orbits. Set of all such points has dimension $2$ and will be denoted $\mathcal N$. For example, points $N = (1,0,0,0)^T$ and $N \cdot R$ lie in $\mathcal N$ for all $R$. 
 \item Points with coordinates $(x_1,x_2,x_3,x_4)$ with $x_1^2 + x_2^2 = x_3^2 + x_4^2 = \frac {|x|^2} 2$. Orbits of such points have length $2 k \pi |x|$. Their orbits will be called {\bf equatorial} orbits. Set of all such points has dimension $3$ and will be denoted $\mathcal M$. For example, points $M = (\frac 1 {\sqrt 2},0,\frac 1 {\sqrt 2},0)^T$ and $M \cdot R$ lie in $\mathcal M$ for all $R$. 
 \item Points with coordinates $(x_1,x_2,x_3,x_4)$ not mentioned in previous items. Orbits of such points have length $4 k \pi |x|$. These are points of general position.
\end{enumerate}

Now, points of type $1$ and type $2$ can be limits of sequence of points of type $3$ or general position points, but limit of sequence of points of type $1$ is always of type $1$ and limit of sequence of points of type $2$ is always of type $2$. We represent this with following graph, called degeneracy graph for orbits of group $\gk$.
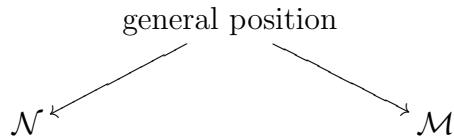
$$
 \xymatrix{
  &\mbox{general position} \ar[ld] \ar[rd]\\
  \mathcal N&&\mathcal M
 }
$$

\begin{lemma}
 \label{global_profile_existence_dim_four_lemma}
 There exist a map $P_G: \reals^4 \to \reals^3$, weight function $\mathfrak w(y)$ and $3\times 3$-matrix function $\hat A(y)$ such that
 \begin{enumerate}
  \item weight function $\mathfrak w$ is $1$-homogeneous and strictly positive;
  \item matrix function $\hat A$ is $0$-homogeneous and uniformly elliptic;
  \item for any $1 \leqslant p,q < \infty$ and any $\mathcal G_1$-invariant function $u$
   \begin{align}
    &\int\limits_{\Omega_R} |u|^q dx = \int\limits_{P_G \Omega_R} |P_G u|^q \mathfrak w(P_G x) d P_G x;\label{global_profile_function_dim_four_estimate}\\
    &\int\limits_{\Omega_R} |\nabla u|^p dx = \int\limits_{P_G \Omega_R} \left ( \hat A(P_G x) \nabla (P_G u),\nabla (P_G u) \right )^{\frac p 2} \mathfrak w(P_G x) d P_G x.\label{global_profile_gradient_dim_four_estimate}
   \end{align}
 \end{enumerate}
 $P_G u$ is called {\bf global profile} of function $u$ and $P_G$ is called a {\bf map providing global profile}.
\end{lemma}

\beginproof
Let $u$ be a $\mathcal G_1$-invariant function. Let 
 $$
  P_0(x_1,x_2,x_3,x_4)=
  \begin{cases}
   (x_1,x_2,x_3,x_4), \mbox{ if } x_1^2 + x_2^2 \leqslant x_3^2 + x_4^2;\\
   (x_3,x_4,x_1,x_2), \mbox{ if } x_1^2 + x_2^2 > x_3^2 + x_4^2.
  \end{cases}
 $$ 
Note that $P_0^{-1}(x)$ consists of exactly two points for any $x \not \in \mathcal M$. As $\mes_4 \mathcal M = 0$, we obtain
\begin{eqnarray*}
 &&\int\limits_{\Omega_R} |u|^q dx = 2 \int\limits_{P_0\Omega_R} |P_0 u|^q d P_0 x\\
 &&\int\limits_{\Omega_R} |\nabla u|^p dx = 2 \int\limits_{P_0\Omega_R} |\nabla(P_0 u)|^p d P_0 x.
\end{eqnarray*}
The map $P_0$ excludes the discrete component of symmetry group $\mathcal G_1$. Namely, $P_0 (\mathcal G_1 x)$ is a connected set for any $x \in \reals^3$. Let $(r_1,\varphi_1)$ be polar coordinates in plane $(x_3,x_4)$. Denote $Q_1(x_1,x_2,x_3,x_4) = (x_1,x_2,r_1,\varphi_1)$ and $P_1 = Q_1 P_0$. Then 
\begin{multline*}
 P_1 \Omega_R = \{ (x_1,x_2,r_1,\varphi_1) : \\
  r_1 > \sqrt{x_1^2 + x_2^2} \geqslant 0;
  (R-1)^2 < r_1^2 + x_1^2 + x_2^2  < (R + 1)^2; 0 < \varphi_1 < 2 \pi \}.
\end{multline*}
and
\begin{multline*}
 P_1 (\mathcal G_1 x^{(0)}) = \{(x_1,x_2,r_1,\varphi_1) : \quad \varphi_1 \in [0,2 \pi], \\ 
 r_1 = r_1(x^{(0)}) = {\rm const}, \quad T_{\varphi_1}^{-1}\vec{ x_1 \\x_2 } = \vec{x_1^{(0)}\\x_2^{(0)}} = {\rm const}\}.
\end{multline*}
 We obtain
\begin{eqnarray*}
 &&\int\limits_{\Omega_R} |u|^q dx = 2 \int\limits_{P_1\Omega_R} |P_1 u|^q r_1 d P_1 x;\\
 &&\int\limits_{\Omega_R} |\nabla u|^p dx = 2 \int\limits_{P_1 \Omega_R} \left[ \left ( \frac {\partial (P_1 u)} {\partial r_1} \right )^2 + \frac 1 {r_1^2}\left ( \frac {\partial (P_1 u)} {\partial \varphi_1} \right )^2 + \sum_{j=1}^{2} \left ( \frac {\partial (P_1 u)} {\partial x_j} \right )^2 \right ]^{\frac p 2} r_1 d P_1 x.
\end{eqnarray*}

Let $\vec{y_1\\y_2} = T_{\varphi_1}^{-1} \vec{x_1\\x_2}$. Define $Q_2(x_1,x_2,r_1,\varphi) = (y_1,y_2,r_1,\varphi_1)$ and $P_2 = Q_2 P_1$. Then 
\begin{multline*}
 P_2 \Omega_R = \{ (y_1,y_2,r_1,\varphi_1) :\\
  r_1 > \sqrt{y_1^2 + y_2^2} \geqslant 0;
  (R-1)^2 < r_1^2 + y_1^2 + y_2^2 < (R + 1)^2; 0 < \varphi_1 < 2 \pi \}
\end{multline*}
and $P_2 (\mathcal G_1 x) = \{r_1={\rm const},y_1={\rm const}, y_2={\rm const}, 0< \varphi_1 < 2 \pi\}$ is a segment of a straight line. Let $Q_3(y_1,y_2,r_1,\varphi_1) = (y_1,y_2,r_1)$ and $P_3 = Q_3 P_2$, so that $P_3 (\mathcal G_1 x)$ consists of exactly one point. As function $u$ is constant with respect to $\varphi_1$, we obtain
\begin{align}
 \int\limits_{\Omega_R} |u|^q dx &= 4 \pi \int\limits_{P_3 \Omega_R} |P_3 u|^q r_1 d P_3(x)\nonumber\\
 \int\limits_{\Omega_R} |\nabla u|^p dx &= 4 \pi \int\limits_{P_3 \Omega_R} \bigg[ \left ( \frac {\partial (P_3 u)} {\partial r_1} \right )^2 + \frac 1 {r_1^2}\left (  (y_{2}\frac {\partial (P_3 u)} {\partial y_{1}} - y_{1}\frac {\partial (P_3 u)} {\partial y_{2}}) \right )^2 + \nonumber \\
 &\qquad\qquad\qquad\qquad+ \left ( \frac {\partial (P_3 u)} {\partial y_1} \right )^2 + \left ( \frac {\partial (P_3 u)} {\partial y_2} \right )^2 \bigg ]^{\frac p 2} r_1 d P_3(x) =\nonumber\\
 &= 2 \pi \int\limits_{P_3\Omega_R} \left( \widehat A (P_3 x) \nabla (P_3 u), \nabla (P_3 u) \right )^{\frac p 2} r_1 d P_3 x, \label{global_profile_gradient_dim_four_eq}
\end{align}
where
\begin{equation}
 \label{matrix_A}
 \widehat A(P_G(x)) = 
 \left ( 
  \begin{array}{ccc}
   1 + \frac {y_2^2} {r_1^2} & - \frac {y_1 y_2} {r_1^2} & 0 \\
   - \frac {y_1 y_2} {r_1^2} & 1 + \frac {y_1^2} {r_1^2} & 0 \\
   0 & 0 & 1  
  \end{array}
 \right ) .
\end{equation}
Let
\begin{equation}
 \label{eigenvectors}
 \xi = \left ( \frac 1 {\sqrt{ y_1^2 + y_2^2}} \vec{y_1 \\ y_2 } \right ) ;
 \qquad
 \eta = \left ( \frac 1 {\sqrt{y_1^2 + y_2^2}} \vec{y_2 \\ -y_1 } \right ) .
\end{equation}
Matrix $\widehat A$ has eigenvalue $1$ of multiplicity $2$ with eigenvectors $(1,0,0)^T, (0, \xi^T)^T$ and eigenvalue $\frac {|P_3 x|^2} {r_1^2}$ of multiplicity $1$ with eigenvector $(0,\eta^T)^T$. Due to inequality $|P_3 x| \geqslant r_1 \geqslant \frac {|P_3 x|} { 2}$ we obtain $|\zeta|^2 \leqslant (\widehat A(P_3 x)\zeta,\zeta) \leqslant 4 |\zeta|^2$ for all $P_3 x \in P_3 \reals^4$, $\zeta \in \reals^{3}$. Taking $P_G = P_3$ and $\mathfrak w(P_G x) = 4 \pi r_1$, we obtain identities \eqref{global_profile_function_dim_four_estimate} and \eqref{global_profile_gradient_dim_four_estimate}.
\endproof

\begin{remark}
 Lemma \ref{global_profile_existence_dim_four_lemma} holds for $u \in \mathcal L_{\mathcal G_0}$ with support bounded away from plane $\{x_1 = x_2 = 0\}$. Indeed, one can put $P_G = Q_1 Q_2 Q_3$, $\mathfrak w(P_G x) = 2 \pi r_1$ and $\widehat A(P_G x)$ given by \eqref{matrix_A}.
\end{remark}

\begin{lemma}
 Let $v \in W^1_p(P_G \Omega_R)$ and $v = 0$ on $P_G(\partial \Omega_R)$. Then there exists a unique function $u \in \mathcal L_{\mathcal G_1}$ such that $P_G u = v$. If $v$ is invariant with respect to the group generated by matrix
$$
 \left (
  \begin{array}{cc}
   T_{\frac {2 \pi} k} & 0 \\
   0 & 1
  \end{array}
 \right ),
$$
then $u \in \lk$. We will say that function $u$ is {\bf recovered from global profile $v$} and denote it $P_G^* v$.
\end{lemma}
\beginproof 
 For any $x \in \reals^4$ the map $Q_0 = P_0$ maps orbit $\mathcal G_1 x$ to a connected set, maps $Q_1$ and $Q_2$ are bijections and $Q_3$ maps $P_2( \mathcal G_1 x)$ to a single point. Let function $w$ be defined by identity $w(r_1,\varphi_1,y_3,y_4) = v(r_1,y_3,y_4)$. Define $u$ by identity
 $$
  u(x_1,x_2,x_3,x_4)=
  \begin{cases}
   (Q^{-1}_{1} Q^{-1}_{2} w)(x_1,x_2,x_3,x_4), \mbox{ if } x_1^2 + x_2^2 \geqslant x_3^2 + x_4^2;\\
   (Q^{-1}_{1} Q^{-1}_{2}w)(x_3,x_4,x_1,x_2), \mbox{ if } x_1^2 + x_2^2 < x_3^2 + x_4^2.
  \end{cases}
 $$ 
 Obviously, function $u \in \mathcal L_{\mathcal G_1}$ and $P_G u = v$.
\endproof

Next corollary is a special case of result of \cite{hebey_vaugon}, but we give it with full proof for reader's convenience.

\begin{corollary}
 \label{embedding_dim_four_lemma}
 Let $1 \leqslant q < p^*_3$. Then $\lk$ is compactly embedded into $L_q(\Omega_R)$.
\end{corollary}

\beginproof
 By \eqref{global_profile_function_dim_four_estimate}, \eqref{global_profile_gradient_dim_four_estimate} $\mathcal L_{\mathcal G_1} \cong W_p^1(P_G \Omega_R, r_1) \cong W_p^1(P_G\Omega_R)$ and $L_q(P_G\Omega_R, r_1) \cong L_q(P_G \Omega_R)$. If $1 \leqslant q<p_3^*$, embedding of $W_p^1(P_G\Omega_R)$ into $L_q(P_G\Omega_R)$ is compact.
\endproof

\begin{corollary}
 \label{Friedrichs_estimate_dim_four_lemma}
 There exists such $c_0$ that for all $u \in \lk$
 \begin{equation}
  \label{lower_estimate_dim_four}
  c_0 R^{\frac 1 p - \frac 1 q} \|u\|_q \leqslant \|\nabla u\|_p, 
 \end{equation}
 where $c_0$ does not depend on $R$.
\end{corollary}

\beginproof
 By \eqref{global_profile_function_dim_four_estimate} and \eqref{global_profile_gradient_dim_four_estimate} we have for any $v_R \in \lk$
 \begin{multline}
  \int\limits_{\Omega_R} |\nabla v_R|^p dx = C \int \limits_{P_G\Omega_R} \left (A(P_G x) \nabla (P_G v_R), \nabla (P_G v_R) \right )^{\frac p 2} r_1 d P_G x \geqslant \\
  \geqslant C R \int \limits_{P_G\Omega_R} |\nabla (P_G v_R)|^p d P_G x \stackrel * \geqslant C R \left ( \int \limits_{P_G\Omega_R} (P_G v_R)^q d P_G x \right )^{\frac p q} \geqslant \\
  \geqslant C R^{1 - \frac p q} \left ( \int \limits_{P_G\Omega_R} (P_G v_R)^q r_1 d P_G x \right )^{\frac p q} = C R^{1 - \frac p q} \left ( \int \limits_{\Omega_R} v_R^q dx \right )^{\frac p q}.
 \end{multline}
 Inequality denoted with * holds by Proposition \ref{Friedrichs_Nazarov_proposition}.
\endproof

\begin{lemma}
 \label{local_profile_dim_four_lemma}
 Let $x_R$ be a sequence of points, such that the type of their orbits does not depend on $R$. Suppose that 

 1) $u_R \in \lgk$ is a sequence of functions supported inside $\bgk(x_R,\rho_R)$, where $\rho_R = R \cdot o_R(1)$ and
 
 2) If $x_R$ is a sequence of points of general position, suppose additionally that ${\rm dist}(x_R, \mathcal N \cup \mathcal M) - \rho_R \geqslant c > 0$.

 Then there exists a sequence of maps $P_R$ such that
 \begin{align}
  \label{local_profile_function_estimate}
  \int\limits_{\Omega_{R}} |u_R|^q dx = {\rm mes}_1 \left (\mathcal G x_R \right ) \int\limits_{P_R \Omega_{R}} |P_R u_R|^q dP_R(x) \cdot (1 + o_R(1)),\\  
  \label{local_profile_gradient_estimate}
  \int\limits_{\Omega_{R}} |\nabla u_R|^p dx = {\rm mes}_1 \left (\mathcal G x_R \right ) \int\limits_{P_R \Omega_{R}} |\nabla (P_R u_R)|^p dP_R(x) \cdot (1 + o_R(1)).
 \end{align}
 The same statement holds if group $\gk$ is replaced by $\mathcal G_0$.

 Function $P_R u_R$ is called a {\bf local profile} of function $u_R$. Such a sequence of maps $P_R$ is called {\bf a sequence of maps providing local profile (SMPLP) centered at $x_R$}.
\end{lemma}

\begin{remark}
 Condition 2 is not automatically implied by condition 1. For example consider sequence $x_R = (\sqrt{R},0,\sqrt{R^2 - R},0)^T$ and $\rho_R = R^{\frac 3 4}$. Then $x_R$ are of general position but ${\rm{dist}} (x_R, \mathcal N) = \sqrt{R} = R^{\frac 3 4} \cdot o_R(1)$.
\end{remark}

\beginproof
First we consider $\mathcal G_0$-invariant functions.

Without loss of generality we can assume that $\left . (x_1^2 + x_2^2) \right |_{x=x_R} \geqslant \left . (x_3^2 + x_4^2) \right |_{x=x_R}$. Thus we can assume that $x_1^2 + x_2^2 \geqslant \frac 1 2(x_3^2 + x_4^2)$ for all $x \in \supp u_R$ if $R$ is large enough. 

We put $\xi_R = \xi |_{P_G x = P_G x_R}$ and $\eta_R = \eta |_{P_G x = P_G x_R}$ with $\xi$ and $\eta$ introduced in \eqref{eigenvectors}.

Freezing coefficients in identity \eqref{global_profile_gradient_dim_four_eq} by Lemma \ref{matrix_freeze_lemma} and Remark \ref{weight_freeze_remark}, we obtain
\begin{equation}
 \label{temp1}
 \begin{split}
  \int\limits_{\Omega_R} |\nabla u_R|^p dx =  2 \pi \int\limits_{P_3\Omega_R} \bigg[ &\left ( \frac {\partial (P_3 u_R)} {\partial r_1} \right )^2 + \left ( \frac 1 {l^2} - 1 \right ) \left ( \sum_{j=1}^{2} (\eta_R)_j \frac {\partial (P_3 u_R)} {\partial y_j} \right )^2 +\\
  & \qquad + \sum_{j=1}^{2} \left ( \frac {\partial (P_3 u_R)} {\partial y_j} \right )^2 \bigg]^{\frac p 2} r_1 d (P_3 x) (1 + o_R(1)),
 \end{split}
\end{equation}
where $l^2 = \left . \frac {r_1^2} {r_1^2 + y_3^2 + y_4^2} \right |_{P_G x = P_G x_R} \in [\frac 1 2, 1]$.

If $x_R \in \mathcal N$, $l=1$ and we can put $P_R = P_G$. 

Suppose that $x_R \not\in \mathcal N$. We put $D = [\xi_R,\eta_R]^{-1}$ and 
$$
 \vec{z_1 \\ z_2} = D \cdot \vec{y_1 \\ y_2}
$$ and let $P_{4,R} x = (z_1, z_2, r_1)$. Note that $P_{4,R} x_R = (z_1, 0, r_1)^T$, and functions $P_{4,R} u_R$ have support in $B(P_{4,R} x, \rho_R)$. We obtain
\begin{multline*}
 \int\limits_{\Omega_R} |\nabla u_R|^p dx = \\
 = 2 \pi \int\limits_{P_{4,R}\Omega_R} \bigg[ \left ( \frac {\partial (P_{4,R} u_R)} {\partial r_1} \right )^2 + \left ( \frac {\partial (P_{4,R} u_R)} {\partial z_1} \right )^2 + \frac 1 {l^2} \left ( \frac {\partial (P_{4,R} u_R)} {\partial z_2} \right )^2 \bigg]^{\frac p 2} r_1 d(P_{4,R} x) (1 + o_R(1)).
\end{multline*}
Let $(r_2,\theta)$ be polar coordinates in plane $(z_1,z_2)$ and let $P_{5,R} x = (r_2, \theta, r_1)$. Note that $P_{5,R} x_R = (r_2,0,r_1)$. We obtain
\begin{multline*}
 \int\limits_{\Omega_R} |\nabla u_R|^p dx = 2 \pi \int\limits_{ P_{5,R} \Omega_R} \Bigg[ \left ( \frac {\partial (P_{5,R} u_R)} {\partial r_1} \right )^2 + \left ( \cos^2\theta + \frac 1 {l^2} \sin^2 \theta \right ) \left ( \frac {\partial (P_{5,R} u_R)} {\partial r_2} \right )^2 + \\
 + 2 \left ( \frac 1 {l^2} - 1 \right ) \frac {\sin\theta\cos\theta} {r_2} \frac {\partial (P_{5,R} u_R)}{\partial r_2} \frac {\partial (P_{5,R} u_R)}{\partial \theta} + \\
 + \left ( \sin^2\theta + \frac 1 {l^2} \cos^2 \theta \right ) \frac 1 {r_2^2} \left ( \frac {\partial (P_{5,R} u_R)} {\partial \theta} \right )^2   \Bigg ]^{\frac p 2} r_2 r_1 d (P_{5,R} x) (1 + o_R(1)).
\end{multline*}
Applying Remark \ref{weight_freeze_remark} and Lemma \ref{matrix_freeze_lemma} once more, we obtain
\begin{equation*}
 \begin{split}
  \int\limits_{\Omega_R} |\nabla u_R|^p dx = 2 \pi \int\limits_{P_{5,R} \Omega_R} \bigg[ & \left ( \frac {\partial (P_{5,R} u_R)} {\partial r_1} \right )^2 + \left ( \frac {\partial (P_{5,R} u_R)} {\partial r_2} \right )^2 + \\
  & \qquad + \frac 1 {l^2} \frac 1 {r_2^2}  \left ( \frac {\partial (P_{5,R} u_R)} {\partial \theta} \right )^2  \bigg]^{\frac p 2} r_2 r_1 d (P_{5,R}x) (1 + o_R(1)).
 \end{split}
\end{equation*}
We set  $\psi = l \theta$ and $P_{6,R} x = (r_2, \psi, r_1)$. Then
\begin{equation*}
 \begin{split}
  \int\limits_{\Omega_R} |\nabla u_R|^p dx = 2 \pi \int\limits_{P_{6,R} \Omega_R} \bigg[ & \left ( \frac {\partial (P_{6,R} u_R)} {\partial r_1} \right )^2 + \left ( \frac {\partial (P_{6,R} u_R)} {\partial r_2} \right )^2 + \\
  & \qquad + \frac 1 {r_2^2}  \left ( \frac {\partial (P_{6,R} u_R)} {\partial \psi} \right )^2 \bigg]^{\frac p 2} r_2 \frac {r_1} l d (P_{6,R} x) (1 + o_R(1)).
 \end{split}
\end{equation*}
Now we unfold polar coordinates $(r_2,\psi)$ to Carthesian coordinates $(Z_1,Z_2)$ and put $P_{7,R} x = (Z_1, Z_2, r_1)$. Freezing weight function $\frac {r_1} l$ with  Lemma \ref{weight_freeze_lemma}, we obtain
\begin{equation*}
 \int\limits_{\Omega_R} |\nabla u_R|^p dx = 2 \pi R \int\limits_{P_{7,R} \Omega_R} \left[ \left ( \frac {\partial (P_{7,R} u_R)} {\partial r_1} \right )^2 + \sum_{j=3}^{4} \left ( \frac {\partial (P_{7,R} u_R)} {\partial Z_j} \right )^2 \right ]^{\frac p 2} d (P_{7,R}x) (1 + o_R(1)).
\end{equation*}
Finally we ensure that supports of profiles of functions with disjunct supports are disjunct by putting $P_{8,R} x = (Y_1, Y_{2},r_1)$, where 
$$
 \vec{Y_1 \\ Y_2} = D^{-1} \vec{Z_1 \\ Z_2}.
$$
Then
\begin{equation}
 \int\limits_{\Omega_R} |\nabla u_R|^p dx = 2 \pi R \int\limits_{P_{8,R} \Omega_R} | \nabla (P_{8,R} u_R)|^p d(P_{8,R} x) (1 + o_R(1)). \label{gradient_profile_estimate_dim_four}
\end{equation}

Similarly we obtain that
\begin{equation}
 \int\limits_{\Omega_R} |u_R|^q dx = 2 \pi R \int\limits_{P_{8,R} \Omega_R} |P_{8,R} u_R|^q d(P_{8,R} x) (1 + o_R(1)). \label{profile_estimate_dim_four}
\end{equation}
This proves that $P_R = P_{8,R}$ is an SMPLP for sequence $x_R$ and group $\mathcal G_0$. 

Now let $u_R$ be a sequence of $\gk$-invariant functions. Let $E_R$ be the connected component of $\supp u_R$ that contains $x_R$, and let $v_R = u_R \cdot \chi_{E_R}$. Note that $v_R$ is $\mathcal G_0$-invariant. If $g x_R \in \mathcal G_0 x_R$ for some $g \in \mathcal H_k$ then $v_R(gx) = v_R(x)$ for all $x \in \supp v_R$. Otherwise, condition 2 provides that supports of $v_R(gx)$ and $v_R(x)$ are disjunct if $g x_R \not\in \mathcal G_0 x_R$. Also, 
\begin{equation}
 \label{k_group_replication}
 u_R = \frac 1 C \sum_{g \in \mathcal H_k} v_R(g x),
\end{equation}
where 
\begin{equation*}
 C = {\#\{g \in \mathcal H_k | g x_R \in \mathcal G_0 x_R\}} = 
\begin{cases}
 k,& \mbox{ if $x_R$ is of general position};\\
 2k,& \mbox{ if } x_R \in \mathcal M;\\
 k^2,& \mbox{ if } x_R \in \mathcal N.
\end{cases}
\end{equation*}

Thus we obtain
\begin{align*}
  \int\limits_{\Omega_{R}} |u_R|^q dx = 2 \pi R \frac {2k^2} C \int\limits_{P_{8,R} \Omega_{R}} |P_{8,R} v_R|^q dP_{8,R} x \cdot (1 + o_R(1)),\\  
  \int\limits_{\Omega_{R}} |\nabla u_R|^p dx = 2 \pi R \frac {2k^2} C \int\limits_{P_{8,R} \Omega_{R}} |\nabla (P_{8,R} v_R)|^p dP_{8,R} x \cdot (1 + o_R(1)),
\end{align*}
where $\frac {2k^2} C$ is the number of connected components of $\gk x_R$. Noting that $\mes_1 \gk x_R = 2 \pi R \frac {2k^2} C$, we obtain \eqref{local_profile_function_estimate} and \eqref{local_profile_gradient_estimate}.
\endproof

\begin{remark}
 Note that SMPLP is not unique. Let $P_R$ be a SMPLP. Suppose that $K \subset \reals^3$ is an open cone that contains $P_R x_R$ and $L_R: K \to \reals ^3$ is a sequence of $C^{1}$ maps with following properties: 1) $L_R$ is $1$-homogenious; 2) $L_R(P_R x_R) = P_R x_R$; 3) $\left . \frac {\partial L_R} {\partial (P_R x)} \right |_{x=x_R} = I_3$. Then $L_R(P_R (\cdot))$ is also SMPLP.  
\end{remark}

We need a way to reconstruct a sequence of functions from a sequence of their local profiles. We have to consider three subcases: 1) general case, i.e. sequence $x_R$ is not equivalent to any sequence lying in $\mathcal M \cup \mathcal N$; 2) $x_R \in \mathcal N$; 3) $x_R \in \mathcal M$.

Define maps $Q_{i,R}, i=4, \ldots, 8$ by identities $Q_{i,R} (P_{i-1,R} x) = P_{i,R}$. These maps are bijections.

Let $U_R$ be a sequence of local profiles centered at sequence $P_R x_R$, i.e $U_R \in \wolp (P_R \Omega_R)$ and $\supp U_R \in B(P_R x_R, \rho_R)$, where $\rho_R = R \cdot o_R(1)$.

In case of general position we suppose additionally that ${\rm dist}(x_R, \mathcal N \cup \mathcal M) - \rho_R \geqslant c > 0$. We can consider a sequence of functions $u_R = P_G^* Q^{-1}_{4,R} \ldots Q^{-1}_{8,R} U_R$. Obviously $P_R u_R = U_R$.

In case of polar orbit sequence we can use the same argument. However, profile must have a nontrivial symmetry group generated by matrix
\begin{equation*}
 \left ( 
 \begin{array}{cc}
  T_{\frac {2 \pi} k} & 0 \\
  0 & 1
 \end{array}
 \right ).
\end{equation*}

In case of equatorial orbit sequence we may still use the same argument, but profile must satisfy a certain nonlinear symmetry. Namely, introduce spherical coordinates $(r,\varphi,\theta)$ in space $(r_1,Y_3,Y_4)$. Then $v_R(r,\varphi,\theta_1) = v_R(r,\varphi,\theta_2)$ must hold if $\theta_1 + \theta_2 = \frac \pi 2$. For matter of convenience we will introduce additional maps so that symmetry would be linear. Namely, let $P_{9,R}(x) = (r,\varphi,\theta)$. Then $P_{9,R} x_R = (r, \varphi, \frac \pi 4)$ and
\begin{gather*}
 \begin{split}
   \int\limits_{\Omega_R} |\nabla u_R|^p dx = 2 \pi k R \int\limits_{P_{9,R} \Omega_R} \left [ \left ( \frac {\partial (P_{9,R} u_R)} {\partial r} \right )^2 + \frac 1 {r^2} \left ( \frac {\partial (P_{9,R} u_R)} {\partial \theta} \right )^2 + \right . \\
   \left . +  \frac 1 {r^2 \sin^2 \theta} \left ( \frac {\partial (P_{9,R} u_R)} {\partial \varphi} \right )^2 \right ]^{\frac p 2} r^2 \sin \theta d (P_{9.R}x) \cdot (1 + o_R(1));
  \end{split} \\
 \int\limits_{\Omega_R} u_R^p dx = 2 \pi k R \int\limits_{P_{9,R} \Omega_R} (P_{9,R} u_R)^q r^2 \sin \theta d (P_{9.R}x) \cdot (1 + o_R(1)).
\end{gather*}
Using Remark \ref{weight_freeze_remark} and Lemmas \ref{weight_freeze_lemma} and \ref{matrix_freeze_lemma} to freeze the coefficients, obtain
\begin{gather*}
 \begin{split} 
  \int\limits_{\Omega_R} |\nabla u_R|^p dx = 2 \pi k R \int\limits_{P_{9,R} \Omega_R} \left [ \left ( \frac {\partial (P_{9,R} u_R)} {\partial r} \right )^2 + \frac 1 {r^2} \left ( \frac {\partial (P_{9,R} u_R)} {\partial \theta} \right )^2 + \right . \\
  \left . + \frac 2 {r^2} \left ( \frac {\partial (P_{9,R} u_R)} {\partial \varphi} \right )^2 \right ]^{\frac p 2} r^2 \frac 1 {\sqrt{2}} d (P_{9.R}x) \cdot (1 + o_R(1));
 \end{split}\\
 \int\limits_{\Omega_R} u_R^p dx = 2 \pi k R \int\limits_{P_{9,R} \Omega_R} (P_{9,R} u_R)^q r^2 \frac 1 {\sqrt{2}} d (P_{9.R}x) \cdot (1 + o_R(1)).
\end{gather*}
Let $\Theta = \theta + \frac \pi 4, \Phi = \frac \varphi {\sqrt{2}}$, unfold spherical coordinates $(r,\Theta,\Phi)$ into carthesian coordinates $(s_1,s_2,s_3)$ and let $P_{10,R}(x) = (s_1,s_2,s_3)$. We have 
\begin{gather*}
 \begin{split}
  \int\limits_{\Omega_R} |\nabla u_R|^p dx = 2 \pi k R \int\limits_{P_{10,R} \Omega_R} \left [ \left ( \frac {\partial (P_{10,R} u_R)} {\partial s_1} \right )^2 + \left ( \frac {\partial (P_{10,R} u_R)} {\partial s_2} \right )^2 + \right . \\
  \left . + \left ( \frac {\partial (P_{10,R} u_R)} {\partial s_3} \right )^2 \right ]^{\frac p 2}  d (P_{10.R}x) \cdot (1 + o_R(1));
 \end{split} \\
 \int\limits_{\Omega_R} u_R^p dx = 2 \pi k R \int\limits_{P_{10,R} \Omega_R} (P_{10,R} u_R)^q d (P_{10.R}x) \cdot (1 + o_R(1)).
\end{gather*}
$P_{10,R}$ is SMPLP while symmetry constraint becomes linear, namely $(P_R u_R)(s_1,s_2,s_3) = \break (P_R u_R)(s_1,s_2,-s_3)$. Now we rotate coordinate system again so that $P_{11,R}x_R = P_{8,R}x_R$ to ensure that profiles of functions with disjunct supports have disjunct supports. Now we can define the recovering map $P_{11,R}^{-1}$ on any sequence of profiles with a certain mirror symmetry.

Note that in all cases a function can be recovered from an axially symmetric profile.

The next lemma is a variant of Proposition \ref{vanishing_proposition} for $\gk$-invariant functions.

\begin{lemma}
 \label{symmetric_vanishing_dim_four_lemma}
 Let $p \in (1, \infty)$ and $q \in (p,p^*_{3})$. Let sequence $u_R \in \lk$ be such that $\|\nabla u_R\|_p \leqslant C R^{\frac 1 p - \frac 1 q}$ and suppose that for some $\rho>0$
 \begin{equation*}
  \sup_{x \in \reals^4} \int\limits_{\mathfrak B_{\gk}(x,\rho)} u_R^q dx = o_R(1).
 \end{equation*}
 Then $\int\limits_{\Omega_R} u_R^q dx = o_R(1)$.
\end{lemma}
\beginproof
 By Lemma \ref{global_profile_existence_dim_four_lemma} we have for $v_R = P_G u_R$ 
 \begin{gather*}  
  \int\limits_{P_G \Omega_R} |\nabla v_R|^p d P_G x \leqslant C_1 R^{ - \frac p q},\\
  \sup_{x \in \Omega_R} \int\limits_{P_G(\mathfrak B_{\gk}(x,\rho))} |v_R|^q d P_G x = R^{-1} \cdot o_R(1).
 \end{gather*}
 Let $w_R = R^{-\frac 1 q} v_R$. Then
 \begin{gather*}  
  \int\limits_{P_G \Omega_R} |\nabla w_R|^p d P_G x \leqslant C_1 ,\\
  \sup_{x \in \Omega_R} \int\limits_{P_G(\mathfrak B_{\gk}(x,\rho))} |w_R|^q d P_G x = o_R(1).
 \end{gather*}
 Sequence $w_R$ satisfies conditions of Proposition \ref{vanishing_proposition}, so
 $$
  \int\limits_{P_G \Omega_R} |w_R|^q d P_G x = o_R(1),
 $$
 and the statement follows by Lemma \ref{global_profile_existence_dim_four_lemma}.
 \endproof

\begin{remark}
 \label{symmetric_nonvanishing_remark_dim_four}
 Suppose that sequence $u_R \in \lgk$ be such that $\|\nabla u_R\|_p \leqslant C R^{\frac 1 p - \frac 1 q}$ and $\int\limits_{\reals^4}|u_R|^q dx \geqslant c >0$ for some $q \in (p,p^*_{3})$ and some $c$ independent of $R$. Then for any $\rho>0$
\begin{equation*}
 \myliminf_{R\to\infty} \sup_{x \in \reals^4}  {\int\limits_{\mathfrak B_{\gk}(x,\rho)} |u_R|^q}  dx >0
\end{equation*}
 Hence sequence $u_R$ has at least one concentration sequence.
\end{remark}

\section{4D case. Separation of concentration sequences}

\begin{lemma}
 \label{cutoff_dim_four_lemma}
 Let $v_R$ be a sequence of $\gk$-invariant functions. Let sequence $x_R$ be a concentration sequence for $v_R$, i.e. for any $\varepsilon>0$ there are such radii $\rho>0$ and $\rho'(R)$ that (\ref{concentration}) holds. Without loss of generality suppose that $\rho'(R)$ satisfies conditions of lemma \ref{local_profile_dim_four_lemma}. Consider $\gk$-invariant cut-off function $\sigma$ such that
\begin{gather*}
 \sigma(x) = 1 \quad \mbox{if}\quad x \in \mathfrak B_{\gk}(x_R, \frac {5 \rho + \rho'(R)} 6)  \\
 \sigma(x) = 1 \quad \mbox{if}\quad x \not\in \mathfrak B_{\gk}(x_R, \frac {\rho + 5 \rho'(R)} 6)\\
 \sigma(x) = 0 \quad \mbox{if}\quad x \in \mathfrak B_{\gk}(x_R,\frac{ 2 \rho'(R) + \rho} 3) \setminus \mathfrak B_{\gk}(x_R, \frac{\rho'(R) + 2 \rho} 3)\\
 |\nabla \sigma| \leqslant \frac {12} {\rho'(R) - \rho }.
\end{gather*}
Then
\begin{gather}
 \int\limits_{\Omega_R} |\sigma v_R|^q dx \geqslant (1 - \varepsilon) \int\limits_{\Omega_R} |v_R|^q dx \label{cutoff_function_estimate} \\ 
 \int\limits_{\Omega_R} |\nabla(\sigma v_R)|^p dx = \int\limits_{\Omega_R} |\nabla v_R|^p dx (1 + o_R(1)). \label{cutoff_gradient_estimate}
\end{gather}

\end{lemma}

\beginproof
Inequality in \eqref{cutoff_function_estimate} holds due to (\ref{concentration}). Further, by the estimate (\ref{vector_inequality}) we have
\begin{multline*}
 \int\limits_{{\Omega_R}} |\nabla (\sigma v_R)|^p \  dx \leqslant \int\limits_{{\Omega_R}} |\nabla v_R|^p \  dx + \\
 + C_1 \left ( \int\limits_{{\Omega_R}} |\nabla v_R| |\nabla \sigma|^{p-1} |v_R|^{p-1} \sigma \  dx + \int\limits_{{\Omega_R}} |v_R|^p |\nabla\sigma|^p \  dx \right ).
\end{multline*}

Estimate the last term:
\begin{multline*}
 \label{srezka_estimate}
 \int\limits_{{\Omega_R}} |v_R|^p |\nabla\sigma|^p \  dx \leqslant C_2 (\rho'(R) - \rho)^{-p} \int\limits_{\supp \nabla\sigma} \!\!\!\! |v_R|^p \  dx \leqslant \\
 \leqslant C_2 (\rho'(R) - \rho)^{-p} (\mes_4 \supp \nabla\sigma)^{1 - \frac p q} \left ( \int\limits_{{\Omega_R}} |v_R|^q \  dx \right )^{\frac p q}.\\
\end{multline*}
Since $\mes_4\supp \nabla \sigma \leqslant C R (\rho'(R))^3$ and $\rho'(R) \to \infty$, we have 
\begin{multline*}
 \int\limits_{{\Omega_R}} |v_R|^p |\nabla\sigma|^p \  dx \leqslant C_3 (\rho'(R))^{(\frac 3 p - \frac 3 q - 1)p} R^{1 - \frac p q} \left ( \int\limits_{{\Omega_R}} |v_R|^q \  dx \right )^{\frac p q} \stackrel*= \\
 \stackrel*= C_3 R^{1 - \frac p q} \left ( \int\limits_{{\Omega_R}} |v_R|^q \  dx \right )^{\frac p q} \cdot o_R(1) \leqslant \int\limits_{\Omega_R} |\nabla v_R|^p dx \cdot o(1). 
\end{multline*}
Relation * holds because  $q < p^*_3$, and last inequality holds due to (\ref{lower_estimate_dim_four}).

By the H\"older inequality the second term is also $\int\limits_{\Omega_R} |\nabla v_R|^p dx \cdot o_R(1)$, and \eqref{cutoff_gradient_estimate} follows.
\endproof

\begin{remark}
 \label{cutoff_dim_four_remark}
 If sequence $v_R$ has two nonequivalent concentration sequences $x_R$ and $y_R$, one can present a cut-off function that separates neighbourhoods of both concentration points from one another and from neighbourhood of infinity. Fix $\varepsilon>0$ and let $\rho_1, \rho'_1(R)$ and $\rho_2, \rho'_2(R)$ be the radii in inequality (\ref{concentration}) for $x_R$ and $y_R$ respectively. Now we use Remark \ref{slow_exterior_radius_remark} to introduce new exterior radii  
$$
 \rho_i''(R) = \min\left \{ \rho_i'(R), \frac 1 4 \min\left \{ x_R - g y_R | g \in \gk \right \} \right \} \leqslant \rho_i'(R), i=1,2.
$$
 Let $\sigma$ and $\tau$ be the cut-off functions constructed with these radii centered at sequences $x_R$ and $y_R$ respectively. Then supports of $1 - \sigma$ and $1 - \tau$ have an empty intersection, and estimates \eqref{cutoff_function_estimate} and \eqref{cutoff_gradient_estimate} are satisfied for cut-off function $\sigma' = \sigma \tau$. 

Cut-off function which separates neighbourhoods of several concentration sequences and neighbourhood of infinity can be constructed similarly.
\end{remark}

\begin{remark}
 Cut-off function $\sigma'$ equal to $\sigma$ on connected components of $\supp\sigma$ that have nonempty intersection with set $A$ and equal to zero elsewhere is called {\bf component of $\sigma$ separating $A$}.
\end{remark}

\begin{lemma}
 \label{cutoff_no_concentration_dim_four_lemma}
 Let $v_R$ be a sequence of $\gk$-invariant functions. Suppose that sequence $x_R$ is not a concentration sequence for $v_R$, i.e. for any $\varepsilon>0$ there are such radii $\rho'(R)$ that (\ref{vanishing_varepsilon}) holds. Without loss of generality suppose that $\rho'(R)$ satisfies conditions of lemma \ref{local_profile_dim_four_lemma}. Consider $\gk$-invariant cut-off function $\sigma$ such that
\begin{gather*}
 \sigma(x) = 0 \quad \mbox{if}\quad x \in \mathfrak B_{\gk}(x_R, \rho'(R))  \\
 \sigma(x) = 1 \quad \mbox{if}\quad x \not\in \mathfrak B_{\gk}(x_R, \frac {2\rho +  \rho'(R)} 3)\\
 |\nabla \sigma| \leqslant \frac {3} {\rho'(R) - \rho }.
\end{gather*}
Then estimates \eqref{cutoff_function_estimate} and \eqref{cutoff_gradient_estimate} hold.
\end{lemma}
Proof of Lemma \ref{cutoff_no_concentration_dim_four_lemma} repeats the proof of Lemma \ref{cutoff_dim_four_lemma}, using inequality \eqref{vanishing_varepsilon} instead of \eqref{concentration}.

\begin{lemma}
 \label{cutoff_profile_lemma}
 Let $\tau$ be a component of $\sigma$ separating $\gk x_R$. Then one can choose external radius $\rho'(R)$ in definition of $\sigma$ such that $\tau u_R$ satisfies conditions of Lemma \ref{local_profile_dim_four_lemma}.
\end{lemma}
\beginproof
 Suppose that condition 1 of Lemma \ref{local_profile_dim_four_lemma} is not satisfied. Then $\rho'(R) \geqslant c R$, and we use remark \ref{slow_exterior_radius_remark} to introduce new exterior radius $\rho''(R) = c \sqrt R$, which satisfies condition 1.

 Suppose now that condition 2 is not satisfied, that is, $x_R$ is of general position, ${\rm dist}(x_R, \mathcal M \cup \mathcal N) \to \infty$ as $R \to \infty$ but the inequality ${\rm dist}(x_R, \mathcal M \cup \mathcal N) - \rho'(R) \geqslant c >0$ does not hold. We use Remark \ref{slow_exterior_radius_remark} to introduce new exterior radius 
 $$
  \rho''(R) = \frac 1 2 {\rm dist}(x_R, \mathcal M \cup \mathcal N)
 $$
 and condition 2 is satisfied.
\endproof

\section{4D case. Construction of solutions}

If we minimize quotient \eqref{J_functional} on $\gk$-invariant functions, the only concentration orbit will be the orbit with minimal length, namely $N \cdot R$ (see \cite{ain2}). We will construct local minimizers of quotient (\ref{J_functional}) concentrating in the neighbourhood of locally (but not globally) minimal orbit. Following the general idea of \cite{byeon}, we require that denominator of quotient \eqref{J_functional} is bounded away from zero in neighbourhood of chosen orbit.

Consider the set
\begin{equation}
 \label{A_varkappa_definition}
 A_\varkappa = \left\{ x \in \reals^n \ : \ {\rm dist}(x,\mathcal M) \leqslant \varkappa |x| \right\},
\end{equation}
where $\varkappa > 0$ is such that $N \cdot R \not\in A_\varkappa$. We minimize $J[u]$ on the set $\wlk$ of functions $v \in \lk$ such that
\begin{gather}
 \int\limits_{\Omega_R} |v|^q dx = 1, \label{WLK_dim_four}\\
 \int\limits_{A_\varkappa \cap \Omega_R} |v|^q dx \geqslant (1 - \delta) \int\limits_{\Omega_R} |v|^q dx \label{delta_inequality_dim_four}
\end{gather} 
where $\delta$ will be chosen later.

\begin{lemma}
 \label{minimum_obtained_dim_four_lemma}
Minimum of functional (\ref{J_functional}) on the set $\wlk$ is obtained.
\end{lemma}

\beginproof
 Since by Lemma \ref{embedding_dim_four_lemma} embedding of $\lk$ into $L_q({\Omega_R})$ is compact, the set $\wlk$ is weakly closed in $\wolp({\Omega_R})$. Functional $J[u]$ on the set $\wlk$ coincides with convex coercive functional $\|\nabla u\|^p_p$. Hence minimum is obtained. 
\endproof

\begin{remark}
 Obviously substitution of $|v|$ into functional (\ref{J_functional}) yields the same value as $v$. So without loss of generality we can assume that $v \geqslant 0$.
\end{remark}

\begin{lemma}
 \label{J_asymp_estimate_dim_four_lemma}
 Let $u_R$ be a sequence of minimizers of functional \eqref{J_functional} on the set $\wlk$. Then
 \begin{equation}
  \label{J_asymp_dim_four}
  J[u_R] \asymp R^{1 - \frac p q}.
 \end{equation}
\end{lemma}
\beginproof
 Lower bound follows from Corollary \ref{Friedrichs_estimate_dim_four_lemma}. 

 Let us obtain the upper bound: let $P_R$ be an SMPLP for sequence $M \cdot R$. Take an arbitrary radially symmetric positive function $V \in \wolp(B_1^3(0))$, consider a sequence of profiles $V_R(P_R x) = V(P_R x - P_R (M \cdot R))$ and recover $\gk$-invariant function $v_R$. This function satisfies inequality (\ref{delta_inequality_dim_four}), and relations \eqref{local_profile_function_estimate} and \eqref{local_profile_gradient_estimate} yield $J[v_R] \leqslant C_1 R^{1 - \frac p q}$.
\endproof

Lemma \ref{J_asymp_estimate_dim_four_lemma} and Remark \ref{symmetric_nonvanishing_remark_dim_four} prove that any sequence of minimizers $u_R$ has concentration sequences. Next three lemmas describe possible concentration sequences.

\begin{lemma}
 (Unique orbit principle)
 \label{unique_orbit_dim_four_lemma}
 Let $u_R$ be a sequence of minimizers of functional \eqref{J_functional} on the set $\wlk$. Then $u_R$ cannot have more that one con\-cen\-tra\-tion sequence inside ${\rm Int} A_\varkappa$, inside ${\rm Ext} A_\varkappa$ and on $\partial A_\varkappa$.
\end{lemma}
\beginproof
 Assume the contrary. For instance, let $x_R$ and $y_R$ be two concentration sequences inside ${\rm Int} A_\varkappa$  with weights $\lambda_1$ and $\lambda_2$ respectively. Remark \ref{cutoff_dim_four_remark} states that there exists a cut-off function $\sigma$ such that estimates \eqref{cutoff_function_estimate} and \eqref{cutoff_gradient_estimate} hold. Let $\sigma_1$ be a component of $\sigma$ separating  $\gk x_R$,  $\sigma_2$ be a component of $\sigma$ separating $\gk y_R$ and $\sigma_0 = \sigma - \sigma_1 - \sigma_2$. Let $P^{(x)}_R$ and $P^{(y)}_R$ be SMPLPs centered at $x_R$ and $y_R$ respectively. By Lemma \ref{cutoff_profile_lemma} we can consider profiles $P^{(x)}_R(\sigma_1 u_R)$ and $P^{(y)}_R(\sigma_2 u_R)$. We can assume that profiles $P^{(x)}_R(\sigma_1 u_R)$ and $P^{(y)}_R(\sigma_2 u_R)$ are defined on the same spherical layer $\omega_R \subset \reals^3$ and their supports are separated from one another. Introduce $b=P^{(x)}_R(\sigma_1 u_R)$ and $c=P^{(y)}_R(\sigma_2 u_R)$ and assume without loss of generality that inequality \eqref{hump_positive_derivative_condition} holds. Applying Lemma \ref{hump_destruction_lemma} to functions $a \equiv 0$, $b$ and $c$, we obtain function $V_R$ defined on $\omega_R$, for which
 \begin{gather*}
  \int\limits_{\omega_R}  |V_R(\xi)|^q d\xi = \int\limits_{\omega_R} (P^{(x)}_R(\sigma_1 u_R))^q + (P^{(y)}_R(\sigma_2 u_R))^q) d\xi,\\
  \int\limits_{\omega_R} |\nabla V_R(\xi)|^p d\xi < \int\limits_{\omega_R} (|\nabla (P^{(x)}_R(\sigma_1 u_R))|^p + |\nabla (P^{(y)}_R(\sigma_2 u_R))|^p) d\xi - \mu.
 \end{gather*}
 We denote by $\widetilde V_R$ the spherical symmetrization of $V_R$ (see \cite{Ka}), and recover a  $\gk$-invariant function $v_R$ from $\widetilde V_R$. By \eqref{local_profile_function_estimate} and \eqref{local_profile_gradient_estimate} we obtain 
\begin{align*}
 \int\limits_{\Omega_R}  v_R^q d x &= \int\limits_{\Omega_R} ((\sigma_1 u_R)^q + (\sigma_2 u_R)^q) d x \cdot (1 + o_R(1)),\\
 \int\limits_{\Omega_R} |\nabla v_R|^p d x &= \left (\int\limits_{\Omega_R} (|\nabla (\sigma_1 u_R)|^p + |\nabla (\sigma_2 u_R)|^p) dx - \mu_1 R^{1 - \frac p q} \right ) \cdot (1 + o_R(1)).
\end{align*}
Then for function $U_R = \sigma_0 u_R + v_R$ we have 
\begin{equation*}
 J[U_R] < (J[\sigma u_R] - \mu_1 R^{1 - \frac p q}) \cdot (1 + o_R(1)) \leqslant \frac 1 {1 - \varepsilon} J[u_R] (1 - \mu_2 + o_R(1))
\end{equation*}
with some $\mu_2>0$, which contradicts the assumption that $u_R$ is a minimizer if $\varepsilon$ is sufficiently small.
Cases of sequences inside ${\rm Ext} A_\varkappa$ and on $\partial A_\varkappa$ are proven similarly.
\endproof

\begin{lemma}
 \label{minimal_orbit_principle_dim_four_lemma}
 (Minimal orbit principle)
 Let $u_R$ be a sequence of minimizers of functional \eqref{J_functional} on the set $\wlk$. Then if con\-cen\-tra\-tion sequence $y_R$ inside ${\rm Ext} A_\varkappa$ exists, it lies in $\mathcal N$, and the concentration sequence $x_R$ inside ${\rm Int} A_\varkappa$ lies in $\mathcal M$. 
\end{lemma}
\beginproof
 Let $y_R$ be the concentration sequence inside ${\rm Ext} A_\varkappa$. Assume that $y_R$ is not equivalent to any sequence in $\mathcal N$. By Lemma \ref{unique_orbit_dim_four_lemma}, $N \cdot R$ is not a concentration sequence. Let $\sigma_1$ be a cut-off function separating neighbourhoods of $\gk x_R$ and infinity, and let $\sigma_2$ be a cut-off function clearing the neighbourhood of $\gk (N \cdot R)$ (these functions are given by Lemmas \ref{cutoff_dim_four_lemma} and \ref{cutoff_no_concentration_dim_four_lemma}). Assume without loss of generality that supports of $1-\sigma_1$ and $1-\sigma_2$ are disjunct. Let $\sigma = \sigma_1 \sigma_2$ and let $\tau$ be a component of $\sigma$ separating $\gk x_R$. Put $\lambda_R =  \int\limits_{\Omega_R} (\tau u_R)^q dx$. Without loss of generality we can assume that $\lambda_R \to \lambda > 0$.

 Let $P^{(x)}_R$ and $P^{(N)}_R$ be SMPLPs at $x_R$ and $N \cdot R$ respectively. By Lemma \ref{cutoff_profile_lemma} we can consider profile $P^{(x)}_R (\tau u_R)$. Denote by $f_R$ its spherical symmetrization. Let $S_R$ and $s_R$ be  sequences of rotations in $\reals^4$ and $\reals^3$ respectively such that $S_R x_R = N \cdot R$ and $s_R P^{(x)}_R x_R = P^{(N)}_R (N \cdot R)$. Consider a function 
\begin{equation*}
 v_R(x) = \left(P^{(N)}_R\right)^* \left ( k^{\frac 1 q} f_R(s_R^{-1} P^{(x)}_R (S_R x)) \right ).
\end{equation*}
Sequence $v_R$ is centered at $N \cdot R$. Since $\mes_1 \left ( \gk x_R \right ) = 4 \pi k R$ and $\mes_1 \left ( \gk (N \cdot R) \right ) = 4 \pi  R$, we obtain
\begin{align*}
 \int\limits_{\Omega_R} v_R^q dx &= \int\limits_{\Omega_R} (\tau u_R)^q dx \cdot (1 + o_R(1)),\\
 \int\limits_{\Omega_R} |\nabla v_R|^p dx &\leqslant k^{1 - \frac p q} \int\limits_{\Omega_R} |\nabla(\tau u_R)|^p dx \cdot (1 + o_R(1)).
\end{align*}
Define $U_R = v_R + (\sigma - \tau)u_R$. Then 
\begin{multline*}
 \int\limits_{\Omega_R} U_R^q dx = \int\limits_{\Omega_R} ((\sigma - \tau)u_R)^q dx + \int\limits_{\Omega_R} v_R^q dx = \int\limits_{\Omega_R} ((\sigma - \tau)u_R)^q dx + \int\limits_{\Omega_R} (\tau u_R)^q dx \cdot (1 + o_R(1)) \geqslant \\
 \geqslant (1 - \varepsilon) \int\limits_{\Omega_R} (\sigma u_R)^q dx \geqslant 1 - 2 \varepsilon,
\end{multline*}
\begin{multline*}
 \int\limits_{\Omega_R} |\nabla U_R|^p dx = \int\limits_{\Omega_R} |\nabla ((\sigma - \tau)u_R)|^p dx + \int\limits_{\Omega_R} |\nabla v_R|^p dx \leqslant \\ 
 \leqslant \int\limits_{\Omega_R} |\nabla (\sigma u_R)|^p dx - \left (1 - \left ( \frac 1 k \right )^{1 - \frac p q}  \right )\int\limits_{\Omega_R} |\nabla(\tau u_R)|^p dx.
\end{multline*}

Using Proposition \ref{Friedrichs_Nazarov_proposition}, we obtain
\begin{multline*}
 \frac \lambda 2 \leqslant \lambda_R = \int\limits_{\Omega_R} (\tau u_R)^q dx = 4 \pi k R  \int\limits_{P^{(x)}_R \Omega_R} (P^{(x)}_R(\tau u_R))^q d P^{(x)}_R x \cdot (1 + o_R(1)) \leqslant \\
 \leqslant 4 \pi k R  \left (\int\limits_{P^{(x)}_R \Omega_R} |\nabla (P^{(x)}_R(\tau u_R))|^p d P^{(x)}_R x \right )^{\frac q p} \cdot (1 + o_R(1)) = \\
 = \left ( 4 \pi k R \right)^{1 - \frac q p} \left( \int\limits_{\Omega_R} |\nabla(\tau u_R)|^p dx \right)^{\frac q p} \cdot (1 + o_R(1)),
\end{multline*}
which means that 
\begin{equation*}
 \int\limits_{\Omega_R} |\nabla(\tau u_R)|^p dx \geqslant \left ( \frac \lambda 2\right )^{\frac p q} \left ( 4 \pi k R \right)^{1 - \frac p q} \cdot (1 + o_R(1)) \geqslant C \lambda^{\frac p q} R^{1 - \frac p q}. 
\end{equation*}
Then due to lemma \ref{J_asymp_estimate_dim_four_lemma} we have
\begin{equation*}
 J[U_R] \leqslant \frac 1 {(1 - 2 \varepsilon)^{\frac p q}} J[u_R] (1 - \mu + o_R(1))
\end{equation*}
with some $\mu > 0$ independent of $\varepsilon$, which contradicts minimality of $u_R$ if $\varepsilon$ is small enough.

The case of ${\rm Int} A_\varkappa$ is considered similarly.
\endproof

\begin{lemma}
 \label{K_negative_dim_four_lemma}
 Let $\delta = \frac{2}{2+k}$, and let $u_R$ be a sequence of minimizers of functional \eqref{J_functional} on $\wlk$. Then a strict inequality holds in (\ref{delta_inequality_dim_four}) if $R$ is large enough.
\end{lemma}
\beginproof
 Assume the contrary: let $u_R$ be such a sequence of minimizers that equality holds in inequality (\ref{delta_inequality_dim_four}). Due to lemma \ref{minimal_orbit_principle_dim_four_lemma} sequence $u_R$ cannot have more than three concentration sequences, and without loss of generality we can assume that sequences inside ${\rm Int} A_\varkappa$ and ${\rm Ext} A_\varkappa$ are $M \cdot R$ and $N \cdot R$ respectively. Denote by $K_R$ the concentration sequence on $\partial A_\varkappa$.

 Lemma \ref{cutoff_dim_four_lemma} provides a cut-off function $\tau^{(N)}$ separating the neighbourhood of orbit of $N \cdot R$ and the neighbourhood of infinity such that
 \begin{eqnarray*}
  &&\int\limits_{\Omega_R} (\tau^{(N)} u_R)^q dx \geqslant 1 - \frac \varepsilon 4 \\
  &&\int\limits_{\Omega_R} |\nabla(\tau^{(N)} u_R)|^p dx = \int\limits_{\Omega_R} |\nabla u_R|^p dx (1 + o_R(1)).\\
 \end{eqnarray*}
 Cut-off functions $\tau^{(M)}$ and $\tau^{(K)}$ are introduced similarly. By Remark \ref{cutoff_dim_four_remark}, a cut-off function $\tau = \tau^{(N)} \cdot \tau^{(M)} \cdot \tau^{(K)}$ satisfies
 \begin{eqnarray*}
  &&\int\limits_{\Omega_R} (\tau u_R)^q dx \geqslant 1 - \varepsilon \\
  &&\int\limits_{\Omega_R} |\nabla(\tau u_R)|^p dx = \int\limits_{\Omega_R} |\nabla u_R|^p dx (1 + o_R(1)).\\
 \end{eqnarray*}
 Denote by $\sigma$ the component of $\tau$ separating $\gk(N \cdot R) \cup \gk(M \cdot R) \cup \gk(K_R)$ and consider the function $(1 - \sigma) u_R$. By Lemma \ref{cutoff_dim_four_lemma} we have
\begin{equation*}
 \begin{split}
  \supp (1 - \sigma) = \reals^n \setminus \Bigg ( \mathfrak B_{\gk}(N \cdot R, \frac{\rho'(R) + 5 \rho} 6) \cup \mathfrak B_{\gk}(M \cdot R, \frac{\rho'(R) + 5 \rho} 6) \cup \\ \cup \mathfrak B_{\gk}(K_R,\frac{\rho'(R) + 5 \rho} 6) \Bigg ).
 \end{split}
\end{equation*}
Concentration sequences $N \cdot R$, $N \cdot R$ and $K_R$ are cleared, thus by Lemma \ref{unique_orbit_dim_four_lemma} the sequence $(1 - \sigma) u_R$ cannot have concentration sequences at all.
Since 
 $$
  \|\nabla((1 - \sigma) u_R)\|_p \leqslant \|\nabla u_R\|_p + \|\nabla (\sigma u_R)\|_p \leqslant C R^{\frac 1 p - \frac 1 q},
 $$
 function $(1 - \sigma) u_R$ satisfies conditions of Lemma \ref{symmetric_vanishing_dim_four_lemma}. Therefore $\int\limits_{\Omega_R} ((1 - \sigma) u_R)^q dx \to 0$. Thus $\int\limits_{\Omega_R} (\sigma u_R)^q dx = 1 + o_R(1)$ and
 \begin{equation}
  \label{J_cutoff_asymptotic_dim_four}
  J[\sigma u_R] = J[u_R] (1 + o_R(1)).
 \end{equation}

We denote by $\hat\tau^{(N)}$ the component of $\sigma$ separating the neighbourhood of $\gk (N \cdot R)$ and set $w^{(N)}_R = \hat\tau^{(N)} u_R$. We define similarly functions $w^{(M)}$ and $w^{(K)}$ and introduce the "normalized weights":
\begin{align*}
 &a_R = \frac 1 {4 \pi} \int\limits_{\Omega_R} |\nabla w^{(N)}_R|^p dx;  & \alpha_R &= \frac 1 {4 \pi} \int\limits_{\Omega_R}  {w^{(N)}_R}^q dx,\\
 &b_R = \frac 1 {2 k \pi} \int\limits_{\Omega_R} |\nabla w^{(M)}_R|^p dx;  & \beta_R &= \frac 1 {2 k \pi} \int\limits_{\Omega_R}  {w^{(M)}_R}^q dx,\\
 &c_R = \frac 1 {4 k \pi} \int\limits_{\Omega_R} |\nabla w^{(K)}_R|^p dx;  & \gamma_R &= \frac 1 {4 k \pi} \int\limits_{\Omega_R} {w^{(K)}_R}^q dx.
\end{align*}
 Let $P_R^{(N)}$ be an SMPLP\footnote{By Lemma \ref{cutoff_profile_lemma} we can assume that $w^{(N)}_R$ satisfies conditions of Lemma \ref{local_profile_dim_four_lemma}, therefore the profile $P_R^{(N)} (w^{(N)}_R)$ exists.} at $N \cdot R$. We denote by $W^{(N)}_R$ the spherical symmetrization of $P_R^{(N)} (\tau^{(N)} u_R)$ with center at $P_R^{(N)}(N \cdot R)$. Introduce similarly the spherical symmetrizations $W^{(M)}$, $W^{(K)}$ of profiles of $w_R^{(M)}$ and $w_R^{(K)}$ with centers at $P_R^{(M)}(M \cdot R)$ and $P_R^{(K)}(K_R)$ respectively.
By construction we obtain
\begin{gather}
 1 - \varepsilon_R \leqslant 4 \pi \cdot \alpha_R + 2 k \pi \cdot \beta_R + 4 k \pi \cdot \gamma_R \leqslant 1\nonumber\\
 4 \pi \cdot \alpha_R \leqslant \delta\label{OR_constraints}\\
 4 \pi \cdot \alpha_R + 4 k \pi \cdot \gamma_R \geqslant \delta - \varepsilon_R,\nonumber
\end{gather}
where $\varepsilon_R = 1 - \int\limits_{\Omega_R} (\sigma u_R)^q dx \to 0$ as $R \to \infty$. As in proof of Lemma \ref{unique_orbit_dim_four_lemma}, we assume that $W_R^{(N)}$, $W_R^{(M)}$ and $W_R^{(K)}$ are defined on the same spherical layer in $\reals^3$, which we denote by $\omega_R$, and are separated from one another.

Due to estimates \eqref{local_profile_function_estimate} and \eqref{local_profile_gradient_estimate} we have

\begin{align*}
 &\int\limits_{\omega_R} |\nabla W^{(N)}_R|^p d\xi \leqslant \frac { a_R} {R} (1 + o_R(1));  
 & \int\limits_{\omega_R}  (W^{(N)}_R)^q d\xi = \frac { \alpha_R} {R} (1 + o_R(1)),\\
 &\int\limits_{\omega_R} |\nabla W^{(M)}_R|^p d\xi \leqslant \frac { b_R} {R} (1 + o_R(1));
 & \int\limits_{\omega_R}  (W^{(M)}_R)^q d\xi = \frac { \beta_R} {R} (1 + o_R(1)),\\
 &\int\limits_{\omega_R} |\nabla W^{(K)}_R|^p d\xi \leqslant \frac { c_R} {R} (1 + o_R(1));
 & \int\limits_{\omega_R}  (W^{(K)}_R)^q d\xi = \frac { \gamma_R} {R} (1 + o_R(1)),\\
\end{align*}

We show that there exists a function providing a smaller value of functional \eqref{J_functional} than $u_R$.

Sequences $\alpha_R$, $\beta_R$ and $\gamma_R$ cannot tend to zero simultaneously, so, passing to subsequence if necessary, we can assume that at least one of them is bounded away from zero. From sequences $\frac {a_R} {\alpha_R^{\frac p q}}$, $\frac {b_R} {\beta_R^{\frac p q}}$ , $\frac {c_R} {\gamma_R^{\frac p q}}$ we pick ones with denominator bounded away from zero, and from them we select the minimal one.

If it is $\frac {a_R} {\alpha_R^{\frac p q}}$, we denote $m_R = a_R$, $\mu_R = \alpha_R$, $X_R = N \cdot R$, $P_R = P_R{(N)}$ and $V_R =  W^{(N)}_R$. If it is $\frac {b_R} {\beta_R^{\frac p q}}$, we denote $m_R = b_R$, $\mu_R = \beta_R$, $X_R = M \cdot R$, $P_R = P_R{(M)}$ and $V_R =  W^{(M)}_R$. If it is $\frac {c_R} {\gamma_R^{\frac p q}}$, we denote $m_R = c_R$, $\mu_R = \gamma_R$, $X_R = K_R$, $P_R = P_R{(K)}$ and $V_R =  W^{(K)}_R$. Thus $V_R$ is the optimal profile.

Let $s_R^{(N)}$ be a rotation that moves $P_R X_R$ to $P_R^{(N)}N \cdot R$. Introduce a function 
$$
 v_R^{(N)} = \left (\frac {\alpha_R}{\mu_R}\right )^{\frac 1 q} \left ( P_R^{(N)} \right )^* \left (V_R \left ( \left ( s_R^{(N)} \right )^{-1} \xi \right ) \right ).
$$
We introduce function $v_R^{(M)}$ and $v_R^{(K)}$ in a similar way and consider a function
$$
 v_R = v_R^{(N)}+v_R^{(M)}+v_R^{(K)}.
$$
 Layout of "humps"\ of function $v_R$ coincides with layout of "humps"\ of function $\sigma u_R$, masses of these "humps"\ coincide too, and "humps"\ of $v_R$ have an optimal shape. Hence
\begin{gather*}
 \int_{\Omega_R} v_R^q dx = ({4 \pi} \cdot \alpha_R + {2 k \pi} \cdot \beta_R + {4 k \pi} \cdot \gamma_R) (1 + o_R(1));\\
 \begin{split}
 \int_{\Omega_R} |\nabla v_R|^p dx = \left( {4 \pi} \cdot \alpha_R^{\frac p q} + {2 k \pi} \cdot \beta_R^{\frac p q} + {4 k \pi} \cdot \gamma_R^{\frac p q} \right) \frac{m_R}{\mu_R^{\frac p q}} (1 + o_R(1)) \leqslant \\ \leqslant \int_{\Omega_R} |\nabla (\sigma u_R)|^p dx (1 + o_R(1)).
 \end{split}
\end{gather*}
Then, using the first inequality of \eqref{OR_constraints}, we obtain
\begin{equation}
 \label{lower_hump_estimate}  
 J[v_R] = \left( {4 \pi} \cdot \alpha_R^{\frac p q} + {2 k \pi} \cdot \beta_R^{\frac p q} + {4 k \pi} \cdot \gamma_R^{\frac p q} \right) \frac{m_R}{\mu_R^{\frac p q}} (1 + o_R(1))
\end{equation}
and
\begin{equation}
 \label{cutoff_hump_estimate}
 J[v_R] \leqslant J[\sigma u_R] (1 + o_R(1)).
\end{equation}
We set $U_R = v^{(M)}_R$. It has "humps"\ only in the neighbourhood of $M \cdot R$, and
\begin{gather*}
 \int_{\Omega_R} U_R^q dx =  {2 k \pi} \cdot \mu_R (1 + o_R(1)); \quad \int_{\Omega_R} |\nabla U_R|^p dx = {2 k \pi} \cdot m_R (1 + o_R(1)) \\ J[U_R] = \left ( {2 k \pi} \right )^{1 - \frac p q} \frac {m_R}{\mu_R^{\frac p q}} (1 + o_R(1)).
\end{gather*}

Let $\mathfrak M$ be the minimum of function ${4 \pi} \cdot \alpha_R^{\frac p q} + {2 k \pi} \cdot \beta_R^{\frac p q} + {4 k \pi} \cdot \gamma_R^{\frac p q}$ under constraints (\ref{OR_constraints}). Using consecutively \eqref{lower_hump_estimate}, \eqref{cutoff_hump_estimate} and \eqref{J_cutoff_asymptotic_dim_four}, obtain
\begin{multline*}
 J[U_R] \leqslant \frac {\left ( {2 k \pi} \right )^{1 - \frac p q}} {\mathfrak M} \cdot J[v_R] (1 + o_R(1)) \leqslant \\ \leqslant \frac {\left ( {2 k \pi} \right )^{1 - \frac p q}} {\mathfrak M} \cdot J[\sigma u_R] (1 + o_R(1)) \leqslant \frac {\left ( {2 k \pi} \right )^{1 - \frac p q}} {\mathfrak M} \cdot J[u_R] (1 + o_R(1)).
\end{multline*}
Applying Lemma \ref{OR_lemma} with $\mathcal A = {4 \pi}$, $\mathcal B = {2 k \pi}$, $\mathcal C = {4 k \pi}$ and $\delta = \frac 2 {2+k}$, we obtain that $\mathfrak M > \left ( {2 k \pi} \right )^{1 - \frac p q}$ if $\varepsilon_R$ is small enough. Introduce constant $c_R$ such that $\|c_R U_R\|_q^q = 1$. Then $\int_{A_\varkappa} (c_R U_R)^q dx = \int_{\Omega_R} (c_R U_R)^q dx = 1$, and $J[U_R] < J[u_R]$ if $R$ is large enough, which contradicts minimality of $u_R$.
\endproof

\begin{remark}
 \label{one_point_dim_four_remark}
 As a byproduct, we obtain that any sequence $u_R$ has exactly one concentration sequence $x_R$, and  $x_R \in \mathcal M$.
\end{remark}

\begin{lemma}
 \label{solution_dim_four_lemma}
 Let $u_R$ be a minimizer of functional \eqref{J_functional} on $\wlk$. Then $u_R$ is a positive weak solution of boundary value problem \eqref{p_laplacian_eq} if $R$ is large enough.
\end{lemma}

\beginproof
Lemma \ref{K_negative_dim_four_lemma} states that constraint $\int_{A_\varkappa} |u_R|^q dx \geqslant (1 - \delta)\int_{\Omega_R} |u_R|^q dx$ is inactive if $R$ is large enough. Hence $u_R$ is a local minimizer on the set 
\begin{equation*}
 \left \{ u \in \wolp({\Omega_R}) \left | \int\limits_{\Omega_R} |u|^q dx = 1 \right . \right \}.
\end{equation*}

By the Lagrange theorem for some $\lambda \in \reals$ following identity holds:
\begin{equation}
 \label{euler_eq}
 \int\limits_{\Omega_R} \left ( (|\nabla u_R|^{p-2} \nabla u_R, \nabla h) - \lambda |u_R|^{q-2} u_R h  \right ) dx = 0 \quad \forall \quad h \in \lk.
\end{equation}

Due to principle of symmetric criticality (see \cite{palais}) $u_R$ is a nonnegative solution to boundary value problem
\begin{equation}
 \label{euler_equation_w_lagrange_multiplier}
 -\Delta_p u_R = \lambda u_R^{q-1}\quad\mbox{in}\quad{\Omega_R};\qquad u_R=0 \quad\mbox{on}\quad \partial{\Omega_R}.
\end{equation}

Notice that $u_R \geqslant 0$ is a super-$p$-harmonic function, and due to the Harnack inequality for $p$-harmonic functions (\cite{trudinger}) it is positive in ${\Omega_R}$. As left-hand side and right-hand side of equation in (\ref{euler_equation_w_lagrange_multiplier}) have different degrees of homogeneity, a function $v = c u_R$ with some appropriate constant $c$ is a solution to boundary problem (\ref{p_laplacian_eq}).
\endproof

We have thus proven a following
\begin{thm}
 \label{existence_dim_four_thm}
 Let $p \in (1,\infty)$, $q \in (p,p^*_3)$, $k \in \mathbb N$, $k \geqslant 2$. Then there exists $R_0 = R_0(p,q,k)$ such that for any $R>R_0$ boundary problem \eqref{p_laplacian_eq} has a positive weak solution concentrating in neighbourhood of $\gk(M \cdot R)$.
\end{thm}

Obviously, solutions provided by this theorem have different concentration sets and are different if $R$ is large enough. We obtain
\begin{thm}
  \label{multiplicity_dim_four_thm}
  Let $p \in (1,\infty)$, $q \in (p,p^*_3)$. Then for any natural $K$ there exists such $R_0 = R_0(p,q,K)$, that for all $R>R_0$ problem (\ref{p_laplacian_eq}) has at least $K$ nonequivalent positive solutions.
 \end{thm} 

\begin{remark}
 If $p<4$, then $p_3^* > p_4^*$, thus Theorem \ref{multiplicity_dim_four_thm} provides the multiplicity of solutions in supercritical case.
\end{remark}

\section{General case with 4D main subspace. Profiles and \break separation of concentration sequences}

In this section we assume that $n \geqslant 6$. We introduce decomposition $\reals^n = \reals^4 \times \reals^{n-4}$ and consider a group $\mathcal G'_k = \gk \times \mathcal O(n-4)$, where groups $\gk$ are defined in previous section. Subspace $F_1 = \reals^{4} \times \{0\}$ will be called ``body'', subspace $F_2 = \{0\} \times \reals^{n-4}$ will be called ``tail''. 

Structure of orbits here is more complicated. There are five classes of orbits.

Points in $F_1$ are classified as before, namely
\begin{enumerate}
 \item Points that have coordinates $(x_1,x_2,0,0,0,\ldots,0)$ or $(0,0,x_3,x_4,0,\ldots,0)^T$. Orbit of any such point has dimension $1$ and length $4 \pi |x|$. Set of all such points will be denoted $\mathcal N$. For example, points $N = (1,0,0,0,0,\ldots,0)^T$ and $N \cdot R$ lie in $\mathcal N$ for all $R$.
 \item Points with coordinates $(x_1,x_2,x_3,x_4,0,\ldots,0)^T$ with $x_1^2 + x_2^2 = x_3^2 + x_4^2 = \frac {|x|^2} 2$. Orbit of any such point has dimension $1$ and length $2 k \pi |x|$. Set of all such points will be denoted $\mathcal M$. For example, points $M = (\frac 1 {\sqrt 2},0,\frac 1 {\sqrt 2},0,0,\ldots,0)^T$ and $M \cdot R$ lie in $\mathcal M$ for all $R$.
 \item Points with coordinates $(x_1,x_2,x_3,x_4,0,\ldots,0)$ not mentioned in previous items. Orbit of any such point has dimension $1$ and length $4 k \pi |x|$.
\end{enumerate}

Other points are divided into two more classes:
\begin{enumerate}
 \item[4.] Points whose orbits lie only in "tail"{} subspace. These orbits have dimension $n-5 \geqslant 1$. If $n=6$, we note in addition that the length of such orbits equals $2 \pi |x|$. Set of all such points will be denoted $\mathcal N_0$. For example, points $N_0 = (0,0,0,0,1,0,\ldots,0)^T$ and $N_0 \cdot R$ lie in $\mathcal N_0$ for all $R$.
 \item[5.] Points of general position, whose orbits lie neither in "tail"{} nor in "body"{} subspace. Dimension of such orbits is $n-4 \geqslant 2$.
\end{enumerate}
 
The degeneracy graph in this case is following:
$$
 
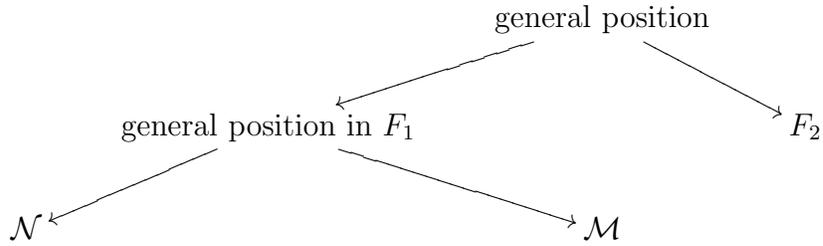
\xymatrix{
  &&\mbox{general position} \ar[ld] \ar[rd]&\\
  &\mbox{general position in }F_1 \ar[ld] \ar[rd]&&F_2\\
  {\mathcal N}&&{\mathcal M}&
 }
$$

First we provide global and local profiles. 

\begin{lemma}
 \label{global_profile_existence_dim_six_and_more_lemma}
 There exist maps $P_{G,1}$, $P_{G,2}$ and $P_{G,3}$,  weight functions $\mathfrak w_1(P_{G,1} x)$, $\mathfrak w_2(P_{G,2} x)$ and $\mathfrak w_3(P_{G,3} x)$ and matrix functions $\hat A_1(P_{G,1} x)$, $\hat A_2(P_{G,2} x)$ and $\hat A_3(P_{G,3} x)$  such that:

 1) weight $\mathfrak w_1$ is $1$-homogeneous, nonnegative and satisfies following inequality: 
  \begin{equation*}
   c_1 \cdot {\rm dist}(x, F_2) \leqslant \mathfrak w_1(P_{G,1}x) \leqslant c_2 \cdot {\rm dist}(x, F_2);   
  \end{equation*}
 2) weight $\mathfrak w_2$ is $(n-5)$-homogeneous, nonnegative and satisfies following inequality:
  \begin{equation*}
   c_3 \cdot \left( {\rm dist}(x, F_1)\right )^{n-5} \leqslant \mathfrak w_2(P_{G,2}x) \leqslant c_4 \cdot \left ( {\rm dist}(x, F_1) \right )^{n-5};
  \end{equation*}
 3) $\mathfrak w_3 = \mathfrak w_1 \mathfrak w_2$;

 4) $(n-1)\times(n-1)$-matrix function $\hat A_1(P_{G,1} x)$, $5\times5$-matrix function $\hat A_2(P_{G,2} x)$ and $4\times4$-matrix function $\hat A_3(P_{G,3} x)$ are $0$-homogeneous and uniformly elliptic; 

 5) following identities hold for all $\mathcal G'_1$-invariant functions $u$ and $i=1,2,3$:
 \begin{align}
  &\int\limits_{\supp u } |u|^q dx = \int\limits_{P_{G,i} (\supp u) } |P_{G,i} u|^q \mathfrak w_i(P_{G,i} x) d P_{G,i} x;\label{global_profile_function_dim_six_and_more}\\
  &\int\limits_{\supp u } |\nabla u|^p dx = \int\limits_{P_{G,i} (\supp u) } \left ( \hat A_i(P_{G,i} x) \nabla (P_{G,i} u),\nabla (P_{G,i} u) \right )^{\frac p 2} \mathfrak w_i(P_{G,i} x) d P_{G,i} x. \label{global_profile_gradient_dim_six_and_more}
 \end{align}
\end{lemma}
 
\beginproof
 Introduce $y_1,y_2,r_1$ as in Lemma \ref{global_profile_existence_dim_four_lemma} and put $P_{G,1} x = (y_1,y_2,r_1,x_5,\ldots,x_n)^T$. Then
 \begin{align}
  \int\limits_{\supp u} |u|^q dx &= 4 \pi \int\limits_{P_{G,1} (\supp u)} |P_{G,1} u|^q r_1 d P_{G,1}(x)\nonumber\\
  \int\limits_{\supp u} |\nabla u|^p dx &= 4 \pi \int\limits_{P_{G,1} (\supp u)} \bigg[ \left ( \frac {\partial (P_{G,1} u)} {\partial r_1} \right )^2 + \frac 1 {r_1^2}\left (  (y_{2}\frac {\partial (P_{G,1} u)} {\partial y_{1}} - y_{1}\frac {\partial (P_{G,1} u)} {\partial y_{2}}) \right )^2 + \nonumber \\
  &+ \left ( \frac {\partial (P_{G,1} u)} {\partial y_1} \right )^2 + \left ( \frac {\partial (P_{G,1} u)} {\partial y_2} \right )^2 + \sum_{j=5}^n \left ( \frac {\partial (P_{G,1} u)} {\partial x_j} \right )^2 \bigg ]^{\frac p 2} r_1 d P_{G,1}(x). 
 \end{align}
 Let 
 $$
  \widehat A_1 = \left (
   \begin{array}{cc}
    \widehat A & 0 \\
    0 & I_{n-4}
   \end{array}
  \right ),
 $$
 where $\widehat A$ is given by \eqref{matrix_A}. Putting $\mathfrak w_1 = 4 \pi r_1$, we obtain \eqref{global_profile_function_dim_six_and_more} and \eqref{global_profile_gradient_dim_six_and_more} for $i=1$.

 Now let $r_0 = \sqrt{x_5^2 + \ldots + x_n^2}$ and put $P_{G,2} x = (x_1,\ldots,x_4,r_0)^T$. Then
 \begin{gather}
  \int\limits_{\supp u} |u|^q dx = \mes_{n-5} (\mathbb S^{n-4}) \int\limits_{P_{G,2} (\supp u)} |P_{G,2} u|^q r_0^{n-5} d P_{G,2}(x)\nonumber\\
  \begin{split}
   \int\limits_{\supp u} |\nabla u|^p dx = \mes_{n-5} (\mathbb S^{n-4}) \int\limits_{P_{G,2} (\supp u)} \bigg[ \sum_{j=1}^4 &\left ( \frac {\partial (P_{G,2} u)} {\partial x_j} \right )^2 +  \\  &+ \left ( \frac {\partial (P_{G,2} u)} {\partial r_0} \right )^2 \bigg ]^{\frac p 2} r_0^{n-5} d P_{G,2}(x).
  \end{split} 
 \end{gather}
 Putting $\mathfrak w_2 = \mes_{n-5} (\mathbb S^{n-4}) r_0^{n-5}$ and $\widehat A_2 = I_5$, we obtain \eqref{global_profile_function_dim_six_and_more} and \eqref{global_profile_gradient_dim_six_and_more} for $i=2$.

 Finally, put $P_{G,3} x = (y_1,y_2,r_1,r_0)$. Let 
 $$
  \widehat A_3 = \left (
   \begin{array}{cc}
    \widehat A & 0 \\
    0 & 1
   \end{array}
  \right ),
 $$
 where $\widehat A$ is given by \eqref{matrix_A}. Putting $\mathfrak w_3 = \mathfrak w_1 \mathfrak w_2$, we obtain \eqref{global_profile_function_dim_six_and_more} and \eqref{global_profile_gradient_dim_six_and_more} for $i=3$.
\endproof

\begin{corollary}
 \label{embedding_dim_six_and_more_lemma}
 Let $1 \leqslant q < p^*_{n-1}$. Then $\lgsk$ is compactly embedded into $L_q(\Omega_R)$.
\end{corollary}
This is a particular case of results of \cite{hebey_vaugon}, but for reader's convenience we give a full proof. 

\medskip

\beginproof
 Define cones
 \begin{gather*}
  {\mathcal K_1} = \{x| \sum_{j=1}^4 x_j^2 \geqslant \sum_{j=5}^n x_j^2\} \\
  {\mathcal K_2} = \{x| \sum_{j=1}^4 x_j^2 \leqslant \sum_{j=5}^n x_j^2\} 
 \end{gather*}
 Apply the map $P_{G,1}$ intoduced in Lemma \ref{global_profile_existence_dim_six_and_more_lemma} to the cone ${\mathcal K_1}$. On $P_{G,1} {\mathcal K_1}$ we have $r_1^2 \geqslant y_1^2 + y_2^2$ and $r_1^2 + y_1^2 + y_2^2 \geqslant x_5^2 + \ldots + x_n^2$. Then $r_1^2 \geqslant \frac 1 4 |x|^2$. On the other hand, $r_1 \leqslant |x|$. Thus $2 \pi (R-1) \leqslant \mathfrak w_1 \leqslant 4 \pi (R+1)$.

 Then
 \begin{equation*}
  P_{G,1} (\mathcal L_{\mathcal G_1'}|_{{\mathcal K_1} \cap \Omega_R}) \subset W_p^1(P_{G,1} ({\mathcal K_1} \cap \Omega_R), \mathfrak w_1) \simeq W_p^1(P_{G,1} ({\mathcal K_1} \cap \Omega_R)).
 \end{equation*}
 Similarly
 \begin{equation*}
  P_{G,1} (L_q({\mathcal K_1} \cap \Omega_R)) = L_q(P_{G,1} ({\mathcal K_1} \cap \Omega_R), \mathfrak w_1) \simeq L_q(P_{G,1} ({\mathcal K_1} \cap \Omega_R)).
 \end{equation*}
 Then $\mathcal L_{\mathcal G_1'}|_{{\mathcal K_1} \cap \Omega_R}$ is compactly embedded into $L_q({\mathcal K_1} \cap \Omega_R)$ if $1 \leqslant q<p_{n-1}^*$.
 
 Similarly, on $P_{G,2} {\mathcal K_2}$ we have $\mes_{n-5} (\mathbb S^{n-4}) 2^{-\frac {n-5} 2} (R-1)^{n-5} \leqslant \mathfrak w_2 \leqslant \mes_{n-5} (\mathbb S^{n-4}) (R+1)^{n-5}$. Repeating the arguments of previous paragraph, we obtain that $\mathcal L_{\mathcal G_1'}|_{{\mathcal K_1} \cap \Omega_R}$ is compactly embedded into $L_q({\mathcal K_2} \cap \Omega_R)$ if $1 \leqslant q < p_5^*$. 

 Thus $\mathcal L_{\mathcal G_1'}$ is compactly embedded into $L_q(\Omega_R)$, if $q \in [1, \max\{p_5^*, p_{n-1}^*\}) = [1, p_{n-1}^*)$. 
\endproof

\begin{corollary}
 \label{Friedrichs_estimate_dim_six_and_more_lemma}
 There exists such $c_0$ that for all $u \in \lgsk$
 \begin{equation}
  \label{lower_estimate_dim_six_and_more}
  c_0 R^{\frac 1 p - \frac 1 q} \|u\|_q \leqslant \|\nabla u\|_p, 
 \end{equation}
 where $c_0$ does not depend on $R$.
\end{corollary}

\beginproof
 Repeating the arguments of Lemma \ref{Friedrichs_estimate_dim_four_lemma}, we obtain
 \begin{gather*}
  c_1 R^{\frac q p - 1} \int\limits_{\Omega_R \cap {\mathcal K_1}} |u|^q dx \leqslant \left ( \int\limits_{\Omega_R \cap {\mathcal K_1}} |\nabla u|^p dx \right )^{\frac q p}; \\
  c_2 R^{(n-5)(\frac q p - 1)} \int\limits_{\Omega_R \cap {\mathcal K_2}} |u|^q dx \leqslant \left ( \int\limits_{\Omega_R \cap {\mathcal K_2}} |\nabla u|^p dx \right )^{\frac q p}.
 \end{gather*}
 As $R > 2$, statement of lemma follows immediately.
\endproof

\begin{lemma}
 \label{local_profile_dim_seven_and_more_lemma}
 Let $x_R \in F_1$ be a sequence of points, such that the type of their orbits does not depend on $R$. Suppose that 

 1) $u_R \in \lgk$ is a sequence of functions supported inside $\bgk(x_R,\rho_R)$, where $\rho_R = R \cdot o_R(1)$ and
 
 2) If $x_R$ is a sequence of points of general position in $F_1$, suppose additionally that \break ${\rm dist}(x_R, \mathcal N \cup \mathcal M) - \rho_R \geqslant c > 0$.
 
 Then such a sequence of maps $P_R$ exists that
 \begin{align}
  \label{local_profile_function_estimate_dim_six_and_more}
  \int\limits_{\Omega_{R}} |u_R|^q dx = {\rm mes}_1 \left ( \gsk x_R \right ) \int\limits_{P_R \Omega_{R}} |P_R u_R|^q dP_R(x) \cdot (1 + o_R(1)),\\  
  \label{local_profile_gradient_estimate_dim_six_and_more}
  \int\limits_{\Omega_{R}} |\nabla u_R|^p dx = {\rm mes}_1 \left ( \gsk x_R \right ) \int\limits_{P_R \Omega_{R}} |\nabla (P_R u_R)|^p dP_R(x) \cdot (1 + o_R(1)).
\end{align}

 If $n=6$, relations \eqref{local_profile_function_estimate_dim_six_and_more} and \eqref{local_profile_gradient_estimate_dim_six_and_more} hold also for a concentration sequence $x_R \in F_2$.

 The same statements holds for functions with symmetry group $\mathcal G'_0$.
\end{lemma}
\beginproof
The first statement can be proven by repeating the steps given in the proof of Lemma \ref{local_profile_dim_four_lemma} to global profiles $P_{G,1} u_R$. 

The second statement can be proven by applying Lemma \ref{weight_freeze_lemma} to identities \eqref{global_profile_function_dim_six_and_more} and \eqref{global_profile_gradient_dim_six_and_more}, since $\widehat A_2 = I_5$.
\endproof

Observing that global profiles $P_{G,1} u_R$ and local profiles $P_R u_R$ are radially symmetric with respect to $F_2$ and repeating the arguments given in Section 2, we can recover a sequence of functions from any sequence of axially symmetric local profiles.

\medskip

We can repeat the proof of Lemma \ref{symmetric_vanishing_dim_four_lemma}, using Lemma \ref{global_profile_existence_dim_six_and_more_lemma} instead of Lemma \ref{global_profile_existence_dim_four_lemma} and obtain
\begin{lemma}
 \label{symmetric_vanishing_dim_six_and_more_lemma}
 Let $p \in (1, \infty)$ and $q \in (p,p^*_{n-1})$. Let sequence $u_R \in \mathcal L_{\mathcal G'_1}$ be such that $\|\nabla u_R\|_p \leqslant C R^{\frac 1 p - \frac 1 q}$ and suppose that for some $\rho>0$
 \begin{equation*}
  \sup_{x \in \reals^n} \int\limits_{\mathfrak B_{\mathcal G'_1}(x,\rho)} u_R^q dx = o_R(1).
 \end{equation*}
 Then $\int\limits_{\Omega_R} u_R^q dx = o_R(1)$.
\end{lemma}

\begin{remark}
 \label{symmetric_nonvanishing_remark_dim_six_and_more}
 Suppose that sequence $u_R \in \mathcal L_{\mathcal G'_1}$ be such that $\|\nabla u_R\|_p \leqslant C R^{\frac 1 p - \frac 1 q}$ and $\int\limits_{\reals^n}|u_R|^q dx \geqslant c >0$ for some $q \in (p,p^*_{n-1})$ and some $c$ independent of $R$. Then for any $\rho>0$
\begin{equation*}
 \myliminf_{R\to\infty} \sup_{x \in \reals^n}  {\int\limits_{\mathfrak B_{\mathcal G'_1}(x,\rho)} |u_R|^q}  dx >0.
\end{equation*}
 Thus, sequence $u_R$ has at least one concentration sequence.
\end{remark}

Next two lemmas are analogues of Lemmas \ref{cutoff_dim_four_lemma} and \ref{cutoff_no_concentration_dim_four_lemma} and allow to separate a neighbourhood of concentration sequence or to clear a neighbourhood of non-concentration sequence.

\begin{lemma}
 \label{cutoff_dim_seven_and_more_lemma}
 Let $v_R$ be a sequence of $\gsk$-invariant functions. Let $x_R \in F_1$ be a concentration sequence for $v_R$, i.e. for any $\varepsilon>0$ there are such radii $\rho>0$ and $\rho'(R)$ that (\ref{concentration}) holds. Without loss of generality suppose that $\rho'(R)$ satisfies conditions of Lemma \ref{local_profile_dim_four_lemma}. Consider $\gsk$-invariant cut-off function $\sigma$ such that
\begin{gather*}
 \sigma(x) = 1 \quad \mbox{if}\quad x \in \mathfrak B_{\gsk}(x_R, \frac {5 \rho + \rho'(R)} 6)  \\
 \sigma(x) = 1 \quad \mbox{if}\quad x \not\in \mathfrak B_{\gsk}(x_R, \frac {\rho + 5 \rho'(R)} 6)\\
 \sigma(x) = 0 \quad \mbox{if}\quad x \in \mathfrak B_{\gsk}(x_R,\frac{ 2 \rho'(R) + \rho} 3) \setminus \mathfrak B_{\gsk}(x_R, \frac{\rho'(R) + 2 \rho} 3)\\
 |\nabla \sigma| \leqslant \frac {12} {\rho'(R) - \rho }.
\end{gather*}
Then
\begin{gather}
 \int\limits_{\Omega_R} |\sigma v_R|^q dx \geqslant (1 - \varepsilon)\int\limits_{\Omega_R} |v_R|^q dx \label{cutoff_function_estimate_dim_six_and_more}\\ 
 \int\limits_{\Omega_R} |\nabla(\sigma v_R)|^p dx = \int\limits_{\Omega_R} |\nabla v_R|^p dx (1 + o_R(1)).\label{cutoff_gradient_estimate_dim_six_and_more}
\end{gather}

If $n=6$, relations \eqref{cutoff_function_estimate_dim_six_and_more} and \eqref{cutoff_gradient_estimate_dim_six_and_more} hold also for concentration sequence $x_R \in F_2$. 
\end{lemma}
\beginproof
 We repeat the proof of Lemma \ref{cutoff_dim_four_lemma}, using inequality \eqref{lower_estimate_dim_six_and_more} instead of \eqref{lower_estimate_dim_four}.
\endproof

\begin{lemma}
 \label{cutoff_no_concentration_dim_six_and_more_lemma}
 Let $v_R$ be a sequence of $\gsk$-invariant functions. Suppose that sequence $x_R \in F_1$ is not a concentration sequence for $v_R$, i.e. for any $\varepsilon>0$ there are such radii $\rho'(R)$ that (\ref{vanishing_varepsilon}) holds. Without loss of generality suppose that $\rho'(R)$ satisfies conditions of lemma \ref{local_profile_dim_four_lemma}. Consider $\gsk$-invariant cut-off function $\sigma$ such that
\begin{gather*}
 \sigma(x) = 0 \quad \mbox{if}\quad x \in \mathfrak B_{\gsk}(x_R, \rho'(R))  \\
 \sigma(x) = 1 \quad \mbox{if}\quad x \not\in \mathfrak B_{\gsk}(x_R, \frac {2\rho + \rho'(R)} 3)\\
 |\nabla \sigma| \leqslant \frac {3} {\rho'(R) - \rho }.
\end{gather*}
Then estimates \eqref{cutoff_function_estimate_dim_six_and_more} and \eqref{cutoff_gradient_estimate_dim_six_and_more} hold.

If $n=6$, relations \eqref{cutoff_function_estimate_dim_six_and_more} and \eqref{cutoff_gradient_estimate_dim_six_and_more} hold also for non-concentration sequence $x_R \in F_2$. 

\end{lemma}

\beginproof
 We again repeat the proof of Lemma \ref{cutoff_dim_four_lemma}, using inequality \eqref{vanishing_varepsilon} instead of \eqref{concentration} and inequality \eqref{lower_estimate_dim_six_and_more} instead of \eqref{lower_estimate_dim_four}.
\endproof

\section{General case with 4D main subspace. Construction of solutions}

 Similarly to Section 4, we minimize the functional $J[u]$ on the set 
 $$
  \wlgsk  = \{ v \in \lgsk : \mbox{ inequalities \eqref{WLK_dim_four}, \eqref{delta_inequality_dim_four} hold} \},
 $$ 
 with $A_\varkappa$ defined in \eqref{A_varkappa_definition}, $\varkappa > 0$ is such that $N_0 \cdot R \not\in A_\varkappa$ and $N \cdot R \not\in A_\varkappa$ and $\delta > 0$ will be chosen later. 
 The existence of minimizer $u_R$ is established by repeating the proof of Lemma \ref{minimum_obtained_dim_four_lemma} verbatim.
 Repeating the proof of Lemma \ref{J_asymp_estimate_dim_four_lemma}, we obtain that sequence of minimizers satisfies estimate \eqref{J_asymp_dim_four}.

\begin{lemma}
 \label{no_high_orbits_dim_seven_and_more_lemma}
 Let $u_R$ be a sequence of minimizers. Let $x_R$ is a concentration sequence for $u_R$. Then
 \begin{enumerate}
  \item If $n \geqslant 7$, $x_R$ is equivalent to a sequence $y_R \in F_1$;
  \item If $n = 6$, $x_R$ is equivalent either to a sequence $y_R \in F_1$ or to a sequence $z_R \in F_2$.
 \end{enumerate}
\end{lemma} 
\beginproof
 Define $d_{R,1} = {\rm dist}(x_R,F_1)$ and $d_{R,2} = {\rm dist}(x_R,F_2)$. Let $d_{R,1} \to \infty$ as $R \to \infty$. By Lemma  \ref{point_concentration_lemma} for some $\rho > 0$ we have
 \begin{equation*}
  0 < \frac \lambda 2 \leqslant \int\limits_{\mathfrak B_{\gsk}(x_R, \rho)} u_R^q dx.
 \end{equation*}

 1) Suppose that $\lim_{R \to \infty} d_{R,2} = \infty$. Apply the map $P_{G,3}$ defined in Lemma \ref{global_profile_existence_dim_six_and_more_lemma}. Then we have for $P_{G,3} x \in P_{G,3} (\mathfrak B_{\gsk}(x_R, \rho))$
 \begin{equation*}
  c_1 d_{R,1} \leqslant r_0 \leqslant c_2 d_{R,1}; \quad c_3 d_{R,2} \leqslant r_1 \leqslant c_4 d_{R,2}.
 \end{equation*}
 Therefore
 \begin{equation*}
  \frac \lambda 2 \leqslant C_5 \int\limits_{P_{G,3}( \mathfrak B_{\gsk}(x_R, \rho))} (P_{G,3} u_R)^q r_0^{n-5} r_1 d P_{G,3} x \leqslant C_6 d_{R,2} d_{R,1}^{n-5} \int\limits_{P_{G,3}( \mathfrak B_{\gsk}(x_R, \rho))} (P_{G,3} u_R)^q d P_{G,3} x.
 \end{equation*}
 By Proposition \ref{Friedrichs_Nazarov_proposition} we have
 \begin{equation*}
  \lambda^{\frac p q} \leqslant C_7 d_{R,2}^\frac p q d_{R,1}^{(n-5)\frac p q} \int\limits_{P_{G,3}( \mathfrak B_{\gsk}(x_R, \rho))} |\nabla (P_{G,3} u_R)|^p d P_{G,2} x .
 \end{equation*}
 Then
 \begin{equation*}
  \lambda^{\frac p q} \leqslant C_8 d_{R,2}^{\frac p q -1} d_{R,1}^{(n-5)(\frac p q - 1)} \int\limits_{ \mathfrak B_{\gsk}(x_R, \rho)} |\nabla u_R|^p dx
 \end{equation*}
 and thus 
 \begin{equation}
  \label{high_orbit_concentration_estimate}
  J[u_R] \geqslant \int\limits_{ \mathfrak B_{\gsk}(x_R, \rho)} |\nabla u_R|^p dx \geqslant C_8 \lambda^{\frac p q} d_{R,2}^{1 - \frac p q} d_{R,1}^{(n-5)(1 - \frac p q)}.
 \end{equation}
 Recall that $d_{R,1}^2 + d_{R,2}^2 = R^2$, so either $d_{R,2} \geqslant \frac R 2$, or $d_{R,1} \geqslant \frac R 2$. In both cases the right-hand side of inequality \eqref{high_orbit_concentration_estimate} grows faster than $R^{1 - \frac p q}$, which contradicts the estimate \eqref{J_asymp_dim_four}. If $n=6$, this concludes the proof.

 2) Suppose that $d_{R,2} \leqslant C$ for all $R$. Apply the map $P_{G,2}$ defined in Lemma \ref{global_profile_existence_dim_six_and_more_lemma}. Then
 we have $c_5 R \leqslant r_0 \asymp d_{R,1} \leqslant c_6 R$ for $P_{G,2} x \in P_{G,2} (\mathfrak B_{\gsk}(x_R, \rho))$ and
 \begin{equation*}
  \frac \lambda 2 \leqslant C_1 \int\limits_{P_{G,2}( \mathfrak B_{\gsk}(x_R, \rho))} (P_{G,2} u_R)^q r_0^{n-5} d P_{G,2} x \leqslant C_2 R^{n-5} \int\limits_{P_{G,2}( \mathfrak B_{\gsk}(x_R, \rho))} (P_{G,2} u_R)^q d P_{G,2} x.
 \end{equation*}
 By Proposition \ref{Friedrichs_Nazarov_proposition} we have
 \begin{equation*}
  \lambda^{\frac p q} \leqslant C_3 R^{(n-5)\frac p q} \int\limits_{P_{G,2}( \mathfrak B_{\gsk}(x_R, \rho))} |\nabla (P_{G,2} u_R)|^p d P_{G,2} x .
 \end{equation*}
 Then
 \begin{equation*}
  \lambda^{\frac p q} \leqslant C_4 R^{(n-5)(\frac p q - 1)} \int\limits_{ \mathfrak B_{\gsk}(x_R, \rho)} |\nabla u_R|^p dx
 \end{equation*}
 and thus 
 \begin{equation*}
  J[u_R] \geqslant \int\limits_{ \mathfrak B_{\gsk}(x_R, \rho)} |\nabla u_R|^p dx \geqslant C_5 \lambda^{\frac p q} R^{(n-5)(1 - \frac p q)},
 \end{equation*}
 which contradicts estimate \eqref{J_asymp_dim_four}, if $n \geqslant 7$. 
\endproof

Let us discuss the multidimensional analogues of Lemmas \ref{unique_orbit_dim_four_lemma} -- \ref{solution_dim_four_lemma}.

If $n \geqslant 7$, statements and proofs of all these lemmas hold with group $\gsk$ replacing $\gk$.
If $n=6$, then Lemma \ref{unique_orbit_dim_four_lemma} also holds.

The following statements are a 6D analogues of Lemma \ref{minimal_orbit_principle_dim_four_lemma} and Lemma \ref{K_negative_dim_four_lemma}:
\begin{lemma}
 \label{minimal_orbit_principle_dim_six_lemma}
 (minimal orbit principle)
 Let $u_R$ be a sequence of minimizers of functional \eqref{J_functional} on the set $\wlgsk$. Then if con\-cen\-tra\-tion sequence $y_R$ inside ${\rm Ext} A_\varkappa$ exists, it lies in $\mathcal N_0$, and the concentration sequence $x_R$ inside ${\rm Int} A_\varkappa$ lies in $\mathcal M$. 
\end{lemma}

\begin{lemma}
 \label{K_negative_dim_six_and_more_lemma}
 Let $\delta = \frac{1}{1+k}$, and let $u_R$ be a sequence of minimizers of functional \eqref{J_functional} on $\wlgsk$. Then a strict inequality holds in (\ref{delta_inequality_dim_four}) if $R$ is large enough.
\end{lemma}
The proofs in case $n=6$ are similar to cases $n=4$ and $n \geqslant 7$ with $N_0$ replacing $N$.

Lemma \ref{solution_dim_four_lemma} in the case $n=6$ also holds. Thus, for $n \geqslant 6$ we obtain the following theorem:

\begin{thm}
 \label{existence_dim_seven_and_more_thm}
 Let $p \in (1,\infty)$, $q \in (p,p^*_{n-1})$, $k \in \mathbb N$. Then there exists $R_0 = R_0(n,p,q,k)$ such that for any $R>R_0$ boundary problem \eqref{p_laplacian_eq} has a positive weak solution concentrating in neighbourhood of $M \cdot R$.
\end{thm}

Obviously, solutions provided by this theorem have different concentration sets and are different if $R$ is large enough. We obtain 

\begin{thm}
  \label{multiplicity_dim_seven_and_more_thm}
  Let $p \in (1,\infty)$, $q \in (p,p^*_{n-1})$. Then for any natural $K$ there exists such $R_0 = R_0(n,p,q,K)$, that for all $R>R_0$ problem (\ref{p_laplacian_eq}) has at least $K$ nonequivalent positive solutions.
 \end{thm}

\section{General case with even-dimensional main subspace}

If $n$ is large enough, we can consider more general symmetry groups, leading to solutions different from those presented in previous sections.

Let $n=2m$ and let group $\wgk$ be generated by $n\times n$-matrices
\begin{eqnarray*}
 \mathcal R_\varphi = \left (
  \begin{array}{cccc}
    T_\varphi & 0 & \ldots & 0 \\
    0 & T_\varphi & \ldots & 0 \\
    \vdots & \vdots & \ddots & \vdots\\
    0 & 0 & \ldots & T_\varphi\\
  \end{array}
 \right ), \qquad
 \mathcal S_k = \left (
  \begin{array}{cccc}
    T_{\frac {2 \pi} k} & 0 & \ldots & 0 \\
    0 & I_2 & \ldots & 0 \\
    \vdots & \vdots & \ddots & \vdots\\
    0 & 0 & \ldots & I_2\\
  \end{array}
 \right ),\\
 \mathcal T_\sigma = \left (
  \begin{array}{cccc}
    \delta_{1,\sigma(1)} I_2 & \delta_{1,\sigma(2)} I_2 & \ldots & \delta_{1,\sigma(m)} I_2 \\
    \delta_{2,\sigma(1)} I_2 & \delta_{2,\sigma(2)} I_2 & \ldots & \delta_{2,\sigma(m)} I_2 \\
    \vdots & \vdots & \ddots & \vdots\\
    \delta_{m,\sigma(1)} I_2 & \delta_{m,\sigma(2)} I_2 & \ldots & \delta_{m,\sigma(m)} I_2 \\
  \end{array}
 \right ).
\end{eqnarray*}

The orbit of any point except the origin under the action of group $\wgk$ has dimension $1$.

We will say that vector $x = (x_1,x_2,\ldots,x_{2m-1},x_{2m})^T$ consists of $m$ {\bf pairs} $(x_{2j-1},x_{2j})^T, j = 1, \ldots, m$. Note that for a given $x$ $\mathcal R_\varphi$ and $\mathcal S_k$ conserve the norm of any pair and $\mathcal T_\sigma$ conserves the set of norms of all pairs.

It is reasonable to classify points by a number of nonzero pairs, or, equi\-val\-en\-t\-ly, by dimension of coordinate hyperspace containing one connected component of orbit of point under the action of group $\wgk$. Denote by $\mathcal F_l$ a set of points with exactly $l$ nonzero pairs. Then we can present a following classification:
\begin{enumerate}
 \item All points in $\mathcal F_1$ have length $2 \pi m |x|$. For the sake of uniformity we put $\mathcal M_1 = \mathcal F_1$.
 \item Points in $\mathcal F_l$, $l \geqslant 2$, are split into two subclasses:
 \begin{enumerate}
  \item Points with equal norms of all their pairs. The length of orbit of any such point is $2 \binom m l k^{l-1} \pi |x|$. Set of all such points will be denoted $\mathcal M_l$.
  \item All other points. The length of orbit of any such point is not less than $2 m \binom m l k^{l-1} \pi |x|$ and not greater than $2 l! \pi k^{l-1} |x|$.
 \end{enumerate}
\end{enumerate}

\begin{remark}
 There exists $k_0 = k_0(m)$ such that if $k \geqslant k_0$, $l \in \{1,\ldots,m-1\}$, $x \in \mathcal F_l$ and $y \in \mathcal F_{l+1}$ then $\mes_1 (\wgk x) < \mes_1 (\wgk y)$. 
\end{remark}

We give a thorough classification of orbits of group $\wgk$ in the case $n=6$.
\begin{enumerate}
 \item Points in $\mathcal F_1 = \mathcal M_1$, for example, points $M_1 = (1,0,0,0,0,0)$ and $M_1 \cdot R$. Orbits of such points have length $6 \pi |x|$. \footnote{$\mathcal F_1 \cap \mathbb S_6$ consists of exactly one orbit $\wgk M_1$.}.
 \item Points in $\mathcal M_2$, for example, points $M_2 = (\frac 1 {\sqrt 2},0,\frac 1 {\sqrt 2},0,0,0)$ and $M_2 \cdot R$. Orbits of such points have length $6 k \pi |x|$.
 \item Points in $\mathcal F_2 \setminus \mathcal M_2$. Orbits of such points have length $12 k \pi |x|$.
 \item Points in $\mathcal M_3$, for example, points $M_3 = (\frac 1 {\sqrt 3},0,\frac 1 {\sqrt 3},0,\frac 1 {\sqrt 3},0)$ and $M_3 \cdot R$. Orbits of such points have length $2 k^2 \pi |x|$.
 \item Points which have exactly three nonzero pairs and only two of their moduli are equal. Orbits of such points have length $6 k^2 \pi |x|$.
 \item All points not mentioned in previous items. Length of their orbits is $12 k^2 \pi |x|$. These are general position points, and their complement on $\mathbb S_6$ has measure zero.
\end{enumerate}

The degeneracy graph in the case $n=6$ is following:
$$
 
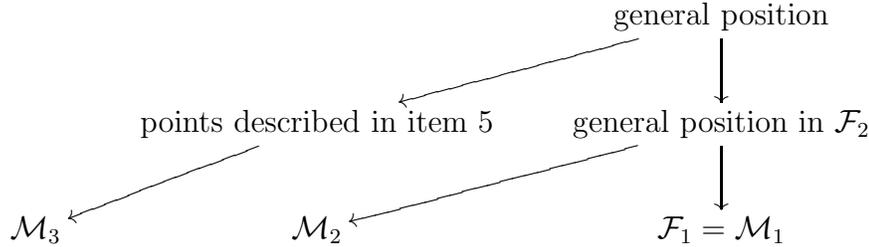
\xymatrix{
  &&\mbox{general position} \ar[ld] \ar[d]&&\\
  &\mbox{points described in item 5} \ar[ld]&\mbox{general position in } \mathcal F_2  \ar[ld] \ar[d]\\
  {\mathcal M_3}&{\mathcal M_2}&{\mathcal F_1 = \mathcal M_1}&
 }
$$

Similarly to the case $n=4$ one can obtain that the following analogues of Lemmas \ref{global_profile_existence_dim_four_lemma} - \ref{local_profile_dim_four_lemma} hold:

\begin{lemma}
 \label{global_profile_existence_dim_even_lemma}
 There exist a the map $P_G: \reals^n \to \reals^{n-1}$, the weight function $\mathfrak w(y)$ and ${(n-1)}\times {(n-1)}$-matrix function $\hat A(y)$ such that
 \begin{enumerate}
  \item the weight function $\mathfrak w$ is $1$-homogeneous and strictly positive;
  \item the matrix function $\hat A$ is $0$-homogeneous and uniformly elliptic;
  \item for any $1 \leqslant p,q < \infty$ and any $\widetilde {\mathcal G_1}$-invariant function $u$ we have
   \begin{align}
    &\int\limits_{\Omega_R} |u|^q dx = \int\limits_{P_G \Omega_R} |P_G u|^q \mathfrak w(P_G x) d P_G x;\label{global_profile_function_dim_even_estimate}\\
    &\int\limits_{\Omega_R} |\nabla u|^p dx = \int\limits_{P_G \Omega_R} \left ( \hat A(P_G x) \nabla (P_G u),\nabla (P_G u) \right )^{\frac p 2} \mathfrak w(P_G x) d P_G x.\label{global_profile_gradient_dim_even_estimate}
   \end{align}
 \end{enumerate}
\end{lemma}

\begin{lemma}
 Let $v \in W^1_p(P_G \Omega_R)$ and $v = 0$ on $P_G(\partial \Omega_R)$. Then there exists a unique function $u \in \mathcal L_{\widetilde {\mathcal G_1}}$ such that $P_G u = v$. If $v$ is invariant with respect to the group generated by block-diagonal matrices
$$
 \left (
  \begin{array}{ccccc}
   T_{\frac {2 \pi} k} & 0 & \ldots & 0 & 0 \\
   0 & I_2 & \ldots & 0 & 0 \\
   \vdots & \vdots & \ddots & \vdots & \vdots \\
   0 & 0 & \ldots & I_2 & 0 \\
   0 & 0 & \ldots & 0 & 1
  \end{array}
 \right ), \ldots,
 \left (
  \begin{array}{ccccc}
   I_2 & 0 & \ldots & 0 & 0 \\
   0 & I_2 & \ldots & 0 & 0 \\
   \vdots & \vdots & \ddots & \vdots & \vdots \\
   0 & 0 & \ldots & T_{\frac {2 \pi} k} & 0 \\
   0 & 0 & \ldots & 0 & 1
  \end{array}
 \right ),
$$
then $u \in \lwgk$.
\end{lemma}

\begin{lemma}
 \label{local_profile_dim_even_lemma}
 Let $x_R$ be a sequence of points, such that the type of their orbits does not depend on $R$. Suppose that 

 1) $u_R \in \lwgk$ is a sequence of functions supported inside $\bwgk(x_R,\rho_R)$, where $\rho_R = R \cdot o_R(1)$ and
 
 2) For $x_R \in \mathcal F_l$, $l \leqslant 2$, suppose additionally that ${\rm dist}(x_R, \mathcal F_1 \cup \ldots \cup \mathcal F_{l-1}) - \rho_R \geqslant c > 0$.

 Then there exists a sequence of maps $P_R$ such that relations \eqref{local_profile_function_estimate} and \eqref{local_profile_gradient_estimate} hold with group $\wgk$ replacing $\gk$.
\end{lemma}

 Fix $l_0 \in \{2,\ldots,m\}$. Consider
 the set
\begin{equation}
 A_\varkappa = \left\{ x \ : \ {\rm dist}(x,\mathcal M_{l_0}) \leqslant \varkappa |x| \right\}.
\end{equation}
where $\varkappa > 0$ is such that $\mathcal F_l \cap A_\varkappa = \emptyset$ for all $l < l_0$. We minimize $J[u]$ on the set $\wlwgk$ of functions $v \in \lwgk$ such that conditions \eqref{WLK_dim_four} and \eqref{delta_inequality_dim_four} hold with $\delta > 0$ to be chosen later.

\begin{proposition}
 \label{A_structure_lemma}
 Let $k\geqslant k_0(m)$. Then orbit of $M_1 \cdot R \in \mathcal M_1$ is the shortest one in ${\rm Ext} A_\varkappa$, orbit of $M_{l_0} \cdot R \in \mathcal M_{l_0}$ is the shortest one in ${\rm Int} A_\varkappa$, and all orbits on $\partial A_\varkappa$ are longer than orbit of $M_{l_0} \cdot R$.
\end{proposition}

Using this proposition, we obtain that Lemma \ref{unique_orbit_dim_four_lemma} holds. Minimal orbit principle also holds in the following form:
\begin{lemma}
 Let $u_R$ be a sequence of minimizers. Then if a con\-cen\-tra\-tion sequence $y_R$ inside ${\rm Ext} A_\varkappa$ exists, it lies in $\mathcal M_1$, if a concentration sequence $z_R$ on $\partial A_\varkappa$ exists, it lies on $\mathcal M_{l_0+1}$, and the concentration sequence $x_R$ inside ${\rm Int} A_\varkappa$ lies in $\mathcal M_{l_0}$.
\end{lemma}

The following lemma is analogue of Lemma \ref{K_negative_dim_four_lemma}:

\begin{lemma}
 Let $\delta = \frac{m}{m+\binom m l k^{l-1}}$, and let $u_R$ be a sequence of minimizers of functional \eqref{J_functional} on $\wlk$. Then a strict inequality holds in (\ref{delta_inequality_dim_four}) if $R$ is large enough.
\end{lemma}

Lemma \ref{solution_dim_four_lemma} and Remark \ref{one_point_dim_four_remark} hold, and thus the following theorem holds:

\begin{thm}
 \label{existence_dim_high_thm}
 Let $p \in (1,\infty)$, $q \in (p,p^*_{n-1})$, $k \in \mathbb N$, $k \geqslant k_0(m)$. Then there exists $R_0 = R_0(p,q,n,k)$ such that for any $R>R_0$ the boundary problem \eqref{p_laplacian_eq} has a positive weak solution concentrating in a neighbourhood of $\wgk(M_{l_0} \cdot R)$.
\end{thm}

Let $n = 2 m + m_0$, $m_0 \geqslant 2$. As we did in Chapter 5, we introduce decomposition $\reals^n = \reals^{2m} \times \reals^{m_0}$ and consider groups ${\widetilde{\mathcal G}}'_k = \wgk \times \mathcal O(m_0)$, $k \in \mathbb N$.  As before, subspace $F_1 = \reals^{2m} \times \{0\}$ will be called ``body'', subspace $F_2 = \{0\} \times \reals^{m_0}$ will be called ``tail''. 

Points in $F_1$ are classified as it was done for $\wgk$, and the degeneracy graph in case $m = 3$ is following:
$$
 
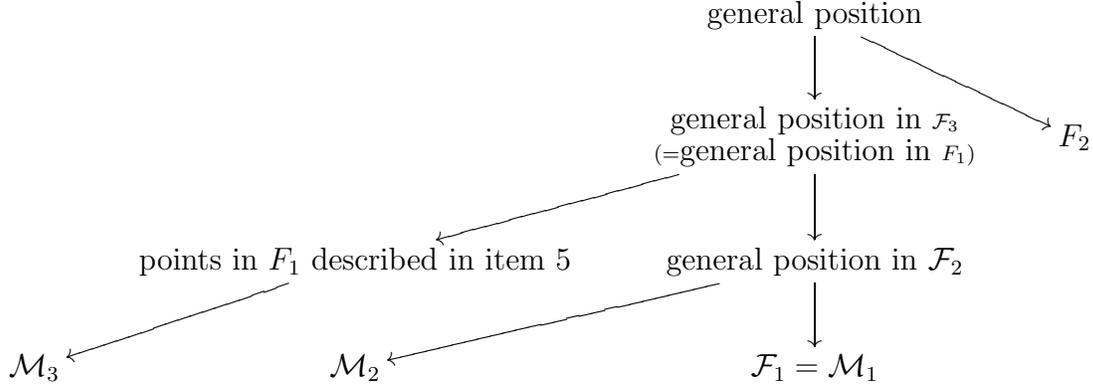
\xymatrix{
  &&\mbox{general position} \ar[d] \ar[rd]&&\\
  &&{\genfrac{}{}{0pt}{}{\mbox{general position in } \mathcal F_3}{(=\mbox{general position in } F_1)}} \ar[ld] \ar[d]&F_2&\\
  &\mbox{points in } F_1 \mbox{ described in item 5} \ar[ld]&\mbox{general position in } \mathcal F_2  \ar[ld] \ar[d]\\
  {\mathcal M_3}&{\mathcal M_2}&{\mathcal F_1 = \mathcal M_1}&
 }
$$

Consider the minimizer of functional \eqref{J_functional} on the set $\wlwgsk$ of functions $v \in \lwgsk$ such that conditions \eqref{WLK_dim_four} and \eqref{delta_inequality_dim_four} hold with $\delta > 0$ to be chosen later. 

Repeating the arguments of chapters 5 and 6, we obtain that the estimate \eqref{J_asymp_dim_four} holds. The following lemma is the analogue of Lemma \ref{no_high_orbits_dim_seven_and_more_lemma}:

\begin{lemma}
 Let $u_R$ be a sequence of minimizers. Let $x_R$ is a concentration sequence for $u_R$. Then
 \begin{enumerate}
  \item For $m_0 \geqslant 3$, $x_R$ is equivalent to a sequence $y_R \in F_1$;
  \item For $m_0 = 2$, $x_R$ is equivalent either to a sequence $y_R \in F_1$ or to a sequence $z_R \in F_2$.
 \end{enumerate}
\end{lemma} 

The case $m_0 \geqslant 3$ can be treated as the case $m_0 = 0$. In case $m_0 = 2$ Lemma \ref{unique_orbit_dim_four_lemma} holds and following analogues of Lemmas \ref{minimal_orbit_principle_dim_four_lemma} and \ref{K_negative_dim_four_lemma} hold:

\begin{lemma}
 (minimal orbit principle)
 Let $u_R$ be a sequence of minimizers of functional \eqref{J_functional} on the set $\wlwgsk$. Then if con\-cen\-tra\-tion sequence $y_R$ inside ${\rm Ext} A_\varkappa$ exists, it lies in $F_2$, and the concentration sequence $x_R$ inside ${\rm Int} A_\varkappa$ lies in $\mathcal M_{l_0}$. 
\end{lemma}

\begin{lemma}
 Let $\delta = \frac{1}{1+\binom m l k^{l-1}}$, and let $u_R$ be a sequence of minimizers of functional \eqref{J_functional} on $\wlk$. Then a strict inequality holds in (\ref{delta_inequality_dim_four}) if $R$ is large enough.
\end{lemma}

Lemma \ref{solution_dim_four_lemma} and Remark \ref{one_point_dim_four_remark} hold, and thus the following theorem holds:

\begin{thm}
 Let $p \in (1,\infty)$, $q \in (p,p^*_{n-1})$, $k \in \mathbb N$, $k \geqslant k_0(m)$. Then there exists $R_0 = R_0(p,q,n,k)$ such that for any $R>R_0$ the boundary problem \eqref{p_laplacian_eq} has a positive weak solution concentrating in a neighbourhood of $\wgsk(M_{l_0} \cdot R)$.
\end{thm}

\section*{Acknowledgements}
Author would like to thank Professor A.I. Nazarov, without whose advices this work would not be possible, and Professor V.G. Osmolovskii, who did the much needed proofreading.

\end{document}